\documentclass{article}
\usepackage{amscd,amsmath,amssymb,amsthm,xspace}
\newcommand{\tiff}{if\textcompwordmark f\xspace}
\newcommand{\theoremname}{\relax}

\newcommand{\inv}{^{-1}}
\newcommand{\ab}{^{ab}}
\newcommand{\abs}[1]{\left|#1\right|}
\newcommand{\Abs}[2][]{{\left\|#2\right\|_{#1}}}
\newcommand{\gint}[1]{\lfloor#1\rfloor}
\newcommand{\sym}[1]{#1\cup#1\inv}
\newcommand{\set}[1]{\{#1\}}
\newcommand{\setst}[2]{\set{\,#1\mid#2\,}}
\renewcommand{\span}[1]{\langle#1\rangle}
\newcommand{\spanst}[2]{\span{\,#1\mid#2\,}}
\newcommand{\into}{\hookrightarrow}
\newcommand{\onto}{\twoheadrightarrow}
\newcommand{\func}[4][\to]{#2\colon#3#1#4}
\newcommand{\lp}[1]{\Z\left[#1,#1\inv\right]}
\newcommand{\ind}[2]{\left[#1:#2\right]}
\newcommand{\sd}{\rtimes}
\newcommand{\ten}{\otimes}
\newcommand{\ds}{\oplus}
\newcommand{\isom}{\cong}
\newcommand{\N}{\mathbb{N}}
\newcommand{\W}{\N\cup\set{0}}
\newcommand{\Z}{\mathbb{Z}}
\newcommand{\Q}{\mathbb{Q}}
\newcommand{\R}{\mathbb{R}}

\newcommand{\II}{\mathcal{I}}
\newcommand{\LL}{\mathcal{L}}
\newcommand{\IM}{I_{max}}
\renewcommand{\Im}{I_{min}}
\newcommand{\LM}{L_{max}}
\newcommand{\Lm}{L_{min}}
\newtheorem{lem}{Lemma}[section]
\newtheorem{prop}[lem]{Proposition}
\newtheorem{cor}[lem]{Corollary}
\newtheorem{thm}[lem]{Theorem}
\newtheorem*{nthm}{\theoremname}
\theoremstyle{definition}
\newtheorem{defin}{Definition}
\theoremstyle{remark}
\newtheorem*{rk}{Remark}
\newtheorem*{notat}{Notation}
\newenvironment{rest}[1]{\renewcommand{\theoremname}{#1}
  \begin{nthm}}{\end{nthm}}

\begin{document}
\title{Strongly $t$-logarithmic $t$-generating sets:\\Geometric
  properties of some soluble groups}
\author{Andrew~D. Warshall\thanks{We thank our advisor, Andrew Casson,
    and Tullia Dymarz for their helpful comments.}\\Yale
  University\\Department of Mathematics\\P.O. Box 208382\\New Haven,
  CT 06520-8382\\\texttt{andrew.warshall@yale.edu}}
\maketitle

\abstract{We introduce the concept of a strongly $t$-logarithmic
  $t$-generating set for a $\lp{t}$-module, which enables us to prove
  that a large class of soluble groups are not almost convex. We also
  prove some results about dead-end depth.}

\section{Introduction}
For an arbitrary finitely generated group $G$ with finite generating
set $A$, the \emph{depth} of an element $g\in G$ is defined to be the
distance (in the word metric with respect to $A$) from $g$ to the
nearest element farther away from the identity. More formally, if
$d(1,g)=n$ then the depth of $g$ is the least integer $d$ such that
$B_g(d)\nsubseteq B_1(n)$. (If there is no such integer, then we say
the depth of $g$ is infinite; this can happen only if $G$ is finite.)

The depth of an element can depend on the choice of generating set; a
classic example of this dependence is $\mathbb{Z}=\span{a}$, which has
depth identically $1$ with respect to the given generating set
$\set{a}$ but in which the depth of $a$ is $2$ with respect to the set
$\set{a^2,a^3}$.  However, for hyperbolic groups it was shown by
Bogopol'ski\u{i} in \cite{B} that the depth is always bounded for any
given generating set.  This is not true for all groups, though; Cleary
and Taback showed in \cite{CT} that the lamplighter group
$\Z_2\wr\Z=\spanst{a,t}{t^2,[t,t^{a^i}],i\in\mathbb{N}}$ has unbounded
depth with respect to the given generating set. In this paper, we give
conditions on a group guaranteeing the existence of generating sets
with unbounded depth, which can be constructed. These conditions also
guarantee that the group is not almost convex with respect to the
generating set so constructed. Using this result, we show that a large
class of soluble groups (including, for example, the soluble
Baumslag-Solitar groups) have unbounded depth with respect to suitable
generating sets, extending a result in \cite{W}.

Let $n\in\N$. A group $G$ is said to be \emph{$n$-almost convex} with
respect to some generating set if for some $N\in\N$ any $g_1$, $g_2\in
G$ at distance at most $n$ from each other are connected by a path in
$G$ of length at most $N$ whose points are at least as close to the
identity as at least one of $g_1$ or $g_2$. If a group is $2$-almost
convex then it is $n$-almost convex for any $n$; this was shown by
Cannon in \cite{C}. In this case we say simply that it is \emph{almost
convex}.  Cannon, Floyd, Grayson and Thurston have shown in
\cite{CFGT} that cocompact lattices in Sol are not almost convex with
respect to any generating set. Later, Miller and Shapiro extended this
result in \cite{MS} to the soluble Baumslag-Solitar groups.  We will
extend these results by showing that the same conditions as in the
preceding paragraph guarantee that the group is not almost convex with
respect to any generating set.

The question has arisen whether the property of having unbounded depth
might be independent of the generating set. In joint work with Riley
(see \cite{RW}), we resolved this question in the negative, producing
a finitely generated group with unbounded depth with respect to one
generating set but depth bounded by $2$ with respect to
another\footnote{The published version of this paper contained an
  error, pointed out by Lehnert; however, the proof of Theorem~3
  contained in Sections~2 and~3 is unaffected.}. In this paper, we
give more examples of groups satisfying the conditions of the
preceding two results but which also have generating sets with respect
to which their depth is bounded. These groups include the lamplighter
group and the soluble Baumslag-Solitar groups.

The paper is organized as follows. In Section~\ref{defs}, we define
some important terms and state the main results. In
Sections~\ref{sunbound} and~\ref{snotac}, we prove that appropriate
conditions on the group imply that it has deep pockets with respect to
some generating set and is not almost convex with respect to any
generating set, respectively. In Section~\ref{sbound}, we prove that a
different set of conditions on the group implies that it lacks deep
pockets with respect to some other generating set. Finally, in
Section~\ref{examples}, we show that a large class of soluble groups
(including the nonabelian soluble Baumslag-Solitar groups) satisfy the
conditions of the first two results. It will follow that those
Baumslag-Solitar groups satisfy the conditions of the third result as
well.

\section{Definitions}\label{defs}
We want to generalize the notion of lamplighter groups. The standard
lamplighter group $L_2=\Z_2\wr\Z$ is given, as relevant here, by the
short exact sequence $0\into\Z_2^\Z\into L_2\onto\Z\onto0$. (By
$\Z_2^\Z$ we mean the collection of biinfinite sequences of elements
of $\Z_2$ such that all but finitely many terms are the identity.)
Since $\Z$ is a free group, this sequence necessarily splits, so $L_2$
may be seen as the semidirect product of $\Z_2^\Z$ with the infinite
cyclic group on one letter, say $t$. In the case of $L_2$, the action
of $t$ on $\Z_2^\Z$ may be taken to be by right shift.

Thus, we may regard $\Z_2^\Z$ as a $\lp{t}$-module, where the action
of $t$ is again by right shift. As a module, it is generated by the
generator of one copy of $\Z_2$; we call that generator $a$. But this
generation is in a stronger sense than the usual, for every element of
$\Z_2^\Z$ may be expressed as $\sum_{i\in\II}t^ia$ for some finite set
$\II$.  (In contrast, the usual notion of a cyclic $\lp{t}$-module
would have every element expressible as $\sum_{i=-\infty}^\infty
n_it^ia$ for all $n_i\in\Z$ and all but finitely many $n_i$ $0$.)

The following definitions generalize the above picture. It is
convenient to include $0$ in all our module generating sets.

\begin{defin}
Let $K$ be an $\lp{t}$-module. Let $0\in A\subseteq K$. Then a formal
expression $\sum_{i=\infty}^\infty t^ia_i$ is a \emph{$t$-word} in $A$
(or just a $t$-word if the choice of $A$ is clear) if all the $a_i\in
A$ and all but finitely many $a_i$ are $0$. More generally, a
\emph{generalized $t$-word} in $A$ is an element of the free abelian
group on $\bigsqcup_{i=-\infty}^\infty t^iA$. Any generalized $t$-word
$v$ (hence, in particular, any $t$-word) \emph{represents} an element
of $K$, denoted $\pi(v)$.
\end{defin}

\begin{rk}
Note that any $t$-word is a generalized $t$-word, but not conversely.
\end{rk}

\begin{defin}
The \emph{length} of a generalized $t$-word is its word length.  The
\emph{length} of a $t$-word is its length as a generalized $t$-word.
A (generalized) $t$-word is \emph{minimal} if it is minimal under
inclusion over all (generalized) $t$-words representing the same
element of $K$.
\end{defin}

\begin{rk}
Note that it is possible, at least \emph{a priori}, for a
(generalized) $t$-word to be minimal yet not of minimal length over
all (generalized) $t$-words representing the same element of $K$.
\end{rk}

\begin{defin}
Let $K$ be a $\lp{t}$-module. Let $0\in A\subseteq K$. The
\emph{$t$-span} of $A$ is the set of all elements of $K$ which are
represented by some $t$-word in $A$.  A subset of $K$ whose $t$-span
is all of $K$ is a \emph{$t$-generating set}. If $A$ is a
$t$-generating set for $K$, it is \emph{symmetrized} if $a\in A$ \tiff
$-a\in A$.
\end{defin}

We will next fix some notations. In the case of the lamplighter group
$\Z_2\wr\Z$, an element $g$ is uniquely specified by two data: its
effect on the position of the lamplighter and its effect on the set of
illuminated lamps. These are generally denoted by $\LL(g)\in\Z$ and
$\II(g)\in P_{fin}(\Z)$, respectively, where $P_{fin}(\Z)$ is the set
of finite subsets of $\Z$. We cannot define a concept equivalent to
$\II(g)$ in general, since many different $t$-words may represent the
same element of $K$. We can, however, make the following definition.

\begin{defin}
Let $K$ be a $\lp{t}$-module and $A$ a symmetrized $t$-generating set
for $K$, with $0\in A$. Let $v$ be a generalized $t$-word in $A$. Then
$\IM(v)$ is the maximal $i$ such that $t^i$ has a nonzero coefficient.
Similarly, $\Im(v)$ is the minimal $i$ such that $t^i$ has a nonzero
coefficient.

Let $k\in K$.  Then $\IM(k)$ is the minimal $\IM(v)$ over all
minimal-length $t$-words $v\in\pi\inv(k)$.  Similarly, $\Im(k)$ is the
maximal $\Im(v)$ over all minimal-length $t$-words $v\in\pi\inv(k)$.
\end{defin}

Note that all of these may in general be $\pm\infty$.

One good property for $A$ to have would be for there to be a unique
$t$-word representing any element of $K$. However, as mentioned above,
this property will hold only rarely. Thus we settle for a weakened
version.

\begin{defin}
Let $K$ be a $\lp{t}$-module and $A$ a symmetrized $t$-generating set
for $K$, with $0\in A$. Then $A$ is \emph{$t$-efficient} if there is
$C\in\N$ with the following property.  Let $\IM$ and $\Im$ be defined
with respect to $A$. Let $v$, $w$ be minimal-length $t$-words in $A$
with $\pi(v)=\pi(w)$.  Then $\abs{\IM(v)-\IM(w)}$,
$\abs{\Im(v)-\Im(w)}<C$.
\end{defin}

\begin{rk}
In the case of the lamplighter group, $\set{a}$ is $t$-efficient with
$C=1$.
\end{rk}

Note that the formal sum of two $t$-words is a generalized $t$-word,
but (even in the case of the lamplighter group) not necessarily a
$t$-word.  The next property we will define for $A$ relates the
possible $t$-word representations of the sum of two elements of $K$ to
the elements' individual $t$-word representations.

\begin{defin}
Let $K$ be a $\lp{t}$-module and $A$ a symmetrized $t$-generating set
for $K$, with $0\in A$.  Then $A$ is \emph{$t$-logarithmic} if there
is $C\in\N$ with the following property.  Let $\IM$ and $\Im$ be
defined with respect to $A$. Then, for all $k_1$, $k_2\in K$,
\[
\IM(k_1+k_2)\le\max(\IM(k_1),\IM(k_2))+C
\]
and
\[
\Im(k_1+k_2)\ge\min(\Im(k_1),\Im(k_2))-C.
\]
\end{defin}

\begin{rk}
In the case of the lamplighter group, $\set{a}$ is $t$-logarithmic
with $C=0$.
\end{rk}

The point of the name is that, with respect to a $t$-logarithmic
$t$-generating set, $\IM(2^ik)\le\IM(k)+Ci$ and
$\Im(2^ik)\ge\Im(k)-Ci$.

The properties of $t$-efficiency and $t$-logarithmicity seem
insufficient for some of our purposes. We will thus define a stronger
property which will give us what we want. To this end, for $n\in\N$
and $w$ a generalized $t$-word, let $\Abs[n]{w}$ denote
$\sum_{i=-\infty}^\infty\max(\Abs{w_i}-n,0)$, where $w_i$ is that
subword of $w$ consisting of letters in $t^iA$.  (Note that
$\Abs[n]{\cdot}$ is not a norm, since it is not multiplicative.) Then
we make the following

\begin{defin}
Let $K$ be a $\lp{t}$-module and $A$ a symmetrized $t$-generating set
for $K$, with $0\in A$.  Then we call $A$ \emph{strongly
  $t$-logarithmic} if for every $m\in\N$ and $n\in\W$ there are
$B_{m,n}$ and $C_{m,n}\in\Z$ with the following property. Let $\IM$
and $\Im$ be defined with respect to $A$. Let $w$ be a nonempty
generalized $t$-word in $A$. Let $w'$ be a minimal $t$-word in $A$
within $n$ of minimal length among all $t$-words representing
$\pi(w)$.  Then
\[
\IM(w')-\IM(w)<B_{m,n}\log(\Abs[m]{w}+1)+C_{m,n}
\]
and
\[
\Im(w)-\Im(w')<B_{m,n}\log(\Abs[m]{w}+1)+C_{m,n}.
\]
\end{defin}

\begin{rk}
In the case of the lamplighter group, $\set{a}$ is strongly
$t$-logarithmic with $B_{m,n}=0$, $C_{m,n}=1$ for all $m$ and $n$.
\end{rk}

\begin{rk}
For the consequences of strong $t$-logarithmicity that we will prove,
the right-hand side of the above inequalities may be replaced by any
function $f(m,n,\Abs[m]{w})$ depending only on $A$ such that, for
every $m$ and $n$, $f$ grows more slowly than every linear function of
$\Abs[m]{w}$.  However, this definition will suffice to include a
significant class of groups.
\end{rk}

We will now justify our assertion that this is a stronger property by
showing that strong $t$-logarithmicity is a generalization of
$t$-ef\-fi\-ciency and $t$-log\-a\-rith\-mic\-ity.

\begin{prop}
If a $t$-generating set is strongly $t$-logarithmic, then it is
$t$-efficient.
\end{prop}

\begin{proof}
The case that $v$ and $w$ are both empty is trivial. Otherwise, in the
definition of strong $t$-logarithmicity, let $w$ also be a
minimal-length $t$-word.  Then $\Abs[1]{w}=0$, so $\IM(w')-\IM(w)$,
$\Im(w)-\Im(w')<C_{1,0}$.
\end{proof}

\begin{prop}
If a $t$-generating set is strongly $t$-logarithmic, then it is
$t$-logarithmic.
\end{prop}

\begin{proof}
The cases that any of $k_1$, $k_2$ or $k_1+k_2$ are $0$ are trivial.
Otherwise, let $w_1\in\pi\inv(k_1)$ and $w_2\in\pi\inv(k_2)$ be
$t$-words such that $\IM(w_1)=\IM(k_1)$ and $\Im(w_2)=\Im(k_2)$.  In
the definition of strong $t$-log\-a\-rith\-mic\-ity, let $w$ be the
formal sum of $w_1$ and $w_2$ and let $w'$ be a minimal-length
$t$-word $\in\pi\inv(k_1+k_2)$.  Note that
\[
\IM(w)\le\max(\IM(w_1),\IM(w_2))=\max(\IM(k_1),\IM(k_2)).
\]
Then $\Abs[2]{w}=0$, so
\begin{multline*}
\IM(k_1+k_2)\le\IM(w')\\<\IM(w)+C_{2,0}\le\max(\IM(k_1),\IM(k_2))+C_{2,0}.
\end{multline*}

The proof for $\Im$ is analogous.
\end{proof}

We will prove the following theorems.

\begin{thm}[deep pockets for some generating set]\label{unbound}
Let $K$ be a nontrivial abelian group and let $T$ be an automorphism
of $K$. Let $G=K\sd\span{t}$, where $t$ acts by $T$, so that $K$ has
the structure of a $\lp{t}$-module. Let $A$ be a finite symmetrized
strongly $t$-logarithmic $t$-generating set for $K$, with $0\in A$.
Let $S=\setst{ta}{a\in A}\cup A$.  Then there are elements of $G$ with
arbitrarily large depth with respect to $S$.
\end{thm}

\begin{thm}[not almost convex]\label{notac}
Let $K$ be a nontrivial abelian group and let $T$ be an automorphism
of $K$. Let $G=K\sd\span{t}$, where $t$ acts by $T$, so that $K$ has
the structure of a $\lp{t}$-module. Suppose there is a finite strongly
$t$-logarithmic $t$-generating set for $K$. Let $S$ be any finite
generating set for $G$.  Then $G$ is not almost convex with respect to
$S$.
\end{thm}

The next theorem requires somewhat different conditions. However,
these conditions are also satisfied by the lamplighter group.

\begin{thm}[no deep pockets for some generating set]\label{bound}
Let $K$ be an abelian group and let $G=K\sd\span{t}$, so that $K$ has
the structure of a $\lp{t}$-module.  Suppose $K$ has a finite
symmetrized $t$-efficient $t$-logarithmic $t$-generating set $A$, with
$0\in A$. Let ${N_A}\ab$ denote the abelianization of the normal
closure of $A$ in $\Z^A\star\span{t}$. Let $\pi$ denote the projection
from ${N_A}\ab$ to $K$.  Suppose there are $I$ and $J\in\N$ such that
for every $n\in\N$ there is $k_n\in K$ such that
\begin{itemize}
\item $\Im(k_n)\ge 1$,
\item $\IM(k_n)\le n-1$ and
\item for any generalized $t$-word $w\in\pi\inv(k_n)$,
  $\abs{N_n(w)}\le I(l(w)-n)+J$, where $l(w)$ is the length of $w$ and
  $N_n(w)$ is the set of $i$ in the range $0<i<n$ such that $w$ does
  not contain any letter of $t^iA$.
\end{itemize}
Let
\[
S=\setst{ta_1t^2a_2t}{a_1,a_2\in\sym{A}\cup\set{0}}\cup A\cup\set{t}.
\]
Then a bound exists on the depth of dead ends in $G$ with respect to
$S$.
\end{thm}

One set of examples of the above theorems (less trivial than the
lamplighter groups) are the soluble Baumslag-Solitar groups $B(1,m)$,
$m>2$.  These may be expressed as $\Z[1/m]\sd\Z$, where the generator
of $\Z$ acts by multiplication by $m$. In Section~\ref{examples}, we
will prove a general result which implies, in particular, that
$A=\set{-\gint{m/2},\dots,\gint{m/2}}$ is a finite strongly
$t$-logarithmic $t$-generating set. Thus in particular it is
$t$-efficient and $t$-logarithmic.

It remains to construct the $k_n$ in the statement of
Theorem~\ref{bound}.  We let
$k_n=\sum_{i=1}^{n-1}m^i=m(m^{n-2}-1)/(m-1)$.  The first two
conditions hold clearly (possibly after changing $n$ by a bounded
amount, which makes no difference). Let $w\in\Z^A$ be a word
representing $k_n$. Regard $w$, read with letters of greatest absolute
value first, as a path through the elements of $\Z[1/m]\subset\R$.
Then, for $0<i<n$, the last multiple of $m^i$ crossed by the path must
be $m^i$ or $2m^i$ modulo $m^{i+1}$. Thus the portion of the path from
the last crossing of a multiple of $m^{i+1}$ (exclusive) to the last
crossing of a multiple of $m^i$ (inclusive) must contain the endpoint
of at least one letter, since by the end of the portion we must be
reading letters of $w$ of absolute value $<m^{i+1}$.  If there is no
letter of $t^iA$ in $w$ then the portion will contain the endpoints of
at least $2$ letters. It follows that $l(w)\ge n-1+N_n(w)$, whence
$N_n(w)\le l(w)-n+1$, so the conditions of Theorem~\ref{bound} are
satisfied for $I=J=1$.

Given $A$, we can define $t$-words without direct reference to
$\lp{t}$-mod\-ules.  Let $F_A$ be the free group on $A$ and let $N_A$
be the normal closure of $A$ in $F_A\star\span{t}$. Then ${N_A}\ab$ is
generated as an abelian group by $t^ia$ for $i\in\Z$ and $a\in A$,
where the action of $t$ on $A$ is by conjugation. We say that an
element of ${N_A}\ab$ is a \emph{$t$-word} if it is of the form
$\sum_{i=-\infty}^\infty t^ia_i$ for $a_i\in A$. Similarly, we can
define generalized $t$-words, $t$-generating sets, and so forth.

\section{Unbounded depth}\label{sunbound}
This section is devoted to the proof of Theorem~\ref{unbound}. Before
we proceed with the proof, we make some more definitions. Let $w\in
F_S$, where $F_S$ is the free group on $S$. Define $\phi(w)$ as
follows.  First, map $w$ to $w'\in F_A\star\span{t}$ in the obvious
way, where $F_A$ denotes the free group on the set $A$.  Express $w'$
as $t^iw''$, with $i\in\Z$ and $w''\in N_A$, where $N_A$, as above,
denotes the normal closure of $A$ in $F_A\star\span{t}$.  Then
${N_A}\ab$ is a free abelian group on $\bigcup_{i=-\infty}^\infty
t^iA$, where $t^iA$ means the image of $A$ under conjugation by $t^i$.
Define $\phi(w)$ to be the image of $w''$ in ${N_A}\ab$.  Let
$\func[\onto]{\sigma}{F_S}{G}$ and $\func[\onto]{\pi}{{N_A}\ab}{K}$ be
the natural projections. Then the following diagram commutes:
\[
\begin{CD}F_S@>\phi>>{N_A}\ab\\@V\sigma_RVV@VV\pi V\\G@>>t^ik\mapsto k>K\end{CD}
\]
where the bottom arrow represents a map taking $t^ik$ to $k$ for
$i\in\Z$ and $k\in K$. (Note that neither horizontal arrow represents
a group homomorphism.)

\begin{proof}[Proof of Theorem~\ref{unbound}]
Let $\func[\onto]{\alpha}{G}{\span{t}}$ be the natural projection.

Let $a\ne0\in A$; such an $a$ exists since $K$ is nontrivial. Choose
$n\in\N$ large and let $g=t^nat^{-2n}at^n\in G$. Then $g\in K$ and
equals $t^na+t^{-n}a$, where we use additive notation since $K$ is
abelian.

I claim that, for $n$ sufficiently large, $t^na+t^{-n}a$ is a minimal
$t$-word.  Note that neither $t^na$ nor $t^{-n}a$ represents $0$,
since $T$ is an automorphism.  Thus they are both minimal-length
$t$-words.  Denote the elements of $K$ they represent by
$\overline{t^na}$ and $\overline{t^{-n}a}$, respectively.  Suppose
$\overline{t^na}=-\overline{t^{-n}a}$.  Then
\[
-n=\IM(t^{-n}a)\ge\IM(\overline{t^{-n}a})=\IM(\overline{t^na})\ge\IM(t^na)-C=n-C
\]
since $A$ is $t$-efficient, where $C$ is (in this paragraph only) as
in the definition of $t$-efficiency.  This is a contradiction for
$n>C/2$, proving the claim.  Furthermore, the length of
$t^na+t^{-n}a$---$2$---is within $1$ of being minimal.

Let $w\in F_S$ be a minimal-length element of $\sigma(w)=g$. Then, by
the above commutative diagram, $\pi(\phi(w))=\overline{t^na+t^{-n}a}$.
Since $A$ is strongly $t$-logarithmic, this implies
\[
n-\IM(\phi(w))<B_{2,1}\log(\Abs[2]{w}+1)+C_{2,1}
\]
and
\[
n+\Im(\phi(w))<B_{2,1}\log(\Abs[2]{w}+1)+C_{2,1},
\]
where $B_{2,1}$ and $C_{2,1}$ are as in the definition of strong
$t$-logarithmicity.  Let $I(w)$ denote the greater of $n-\IM(\phi(w))$
and $n+\Im(\phi(w))$.

We thus get
\[
\Abs[S]{g}\ge\Abs[2]{\phi(w)}+4(n-I(w))>e^{\frac{I(w)-C_{2,1}}{B_{2,1}}}-1+4(n-I(w))>4n-F,
\]
where $F$ is some number dependent only, through $B_{2.1}$, $C_{2,1}$
and $E$, on $A$. (In the above chain of inequalities, the first step
is by the construction of $S$ and $w$'s being of minimal length and
the second step by the preceding paragraph.  The third step is an
application of first-year calculus to the result of the second step,
viewed as a function of $I(w)$.)

Recall that, since $A$ is strongly $t$-logarithmic, it is
$t$-logarithmic.  Let $h\in G$, with $\Abs[S]{h}<(n-E)/(1+C)$, where
$C$ is as in the definition of $t$-logarithmicity and $E$ is the
variable called $C$ in the definition of $t$-efficiency. Let $k\in K$
and $i\in\Z$ be such that $h=t^ik$.  By the construction of $S$, there
is a generalized $t$-word $v''\in\pi\inv(k)$ of length $\le\Abs[S]{h}$
and with $\IM(v'')\le\Abs[S]{h}$ and $\Im(v'')\ge-\Abs[S]{h}$.  Since
$A$ is $t$-logarithmic, it follows that
\[
\IM(k)\le(1+C)\Abs[S]{h}<n-E
\]
and
\[
\Im(k)\ge-(1+C)\Abs[S]{h}>E-n.
\]
Since $A$ is $t$-efficient, there is a (minimal-length) $t$-word
$v\in\pi\inv(k)$ with $\IM(v)\le\IM(k)+E<n$ and
$\Im(v)\ge\Im(k)-E>-n$.  It follows that $hg=t^ik'$, where $k'\in K$
and there is a $t$-word $v'\in\pi\inv(k')$ (namely $v+t^na+t^{-n}a$)
with $\IM(v')=n$ and $\Im(v')=-n$.  Then it is clear that
$\Abs[S]{hg}\le4n<\Abs[S]{g}+F$.  Since $(n-E)/(1+C)$ goes to infinity
as $n$ does, we are done by the Fuzz Lemma from \cite{W}.
\end{proof}

\section{Not almost convexity}\label{snotac}
This section is devoted to proving Theorem~\ref{notac}. This proof is
modeled on that in \cite{MS} for the soluble Baumslag-Solitar groups,
which was in turn modeled on that in \cite{CFGT} for lattices in Sol.
Note that we only use strong $t$-logarithmicity for a restricted class
of words; this is much easier to show, which simplified the work in
those cases, since they could dispense with most of the work in
Section~\ref{examples}.

We begin the proof with the following

\begin{lem}[Triangle Lemma]\label{triangle}
Let $K$ be an $\lp{t}$-module. Let $A$ be a $t$-efficient
$t$-logarithmic $t$-generating set for $K$. Then for every $B\in\R$
there is $D\in\N$ with the following property. Let $w$ be a
generalized $t$-word in $A$. Then there is $n\in\W$ such that $w$ has
more than $Bn$ letters in $\bigcup_{i=\IM(\pi(w))-D-n}^\infty t^iA$
and more than $Bn$ letters in
$\bigcup_{i=-\infty}^{\Im(\pi(w))+D+n}t^iA$.
\end{lem}

\begin{proof}
We will lose nothing if we assume $B\in\N$. Also, we will prove only
the clause involving $\IM$; the proof of the other, involving $\Im$,
is analogous.

Let $C$ be as in the definition of $t$-logarithmicity. Then we choose
subwords $w_1$, $w_2$, \dots of $w$ inductively as follows. Let $w_1$
be a subword formed by choosing up to $BC$ letters of
$\bigcup_{i=\IM(\pi(w))-BC^2-Cj}^\infty t^iA$ for every $j\in\N$.
(Thus, for each $j$ in order, if there are no more than $BC$ letters
in that range not yet chosen then take all of them; otherwise choose
$BC$ of them.)  Then choose $w_2$ similarly as a subword of the
remaining letters, and so on.  Each word $w_i$ is finite since $w$ is
finite.  This process must terminate for the same reason; suppose
$w_k$ is the last nonempty subword. Then $w=\sum_{i=1}^kw_i$, since
each $w_i$ can only be empty if all the letters are already taken.

Thus we have $w=\sum_{i=1}^kw_i$ for some $k\in\W$, where each $w_i$
is a generalized $t$-word. Suppose (for a contradiction) that, for all
$n\in\W$, $w$ has at most $Bn$ letters in
$\bigcup_{i=\IM(\pi(w))-BC^2-n}^\infty t^iA$.  Then, by the
construction of the $w_i$, $\IM(w_i)<\IM(\pi(w))-BC^2-Ci$.  Since each
$w_i$ contains at most $BC$ letters in
$\bigcup_{i=\IM(\pi(w))-BC^2-Cj-C}^{\IM(\pi(w))-BC^2-Cj-1} t^iA$ for
each $j\in\Z$, it is the sum of at most $BC$ $t$-words, which we
denote $w_{i,1}$, \dots, $w_{i,BC}$.  Since each $w_{i,j}$ contains at
most $1$ letter in
$\bigcup_{i=\IM(\pi(w))-BC^2-Cj-C}^{\IM(\pi(w))-BC^2-Cj-1} t^iA$ for
each $j\in\Z$ and satisfies
$\IM(w_{i,j})\le\IM(w_i)<\IM(\pi(w))-BC^2-Ci$, we have
$\IM(\pi(w_{i,j}))<\IM(\pi(w))-BC^2-Ci+C$ by $t$-logarithmicity. Thus
\begin{multline*}
\IM(\pi(w_i))\le\max_{j=1}^{BC}\IM(\pi(w_{i,j}))+BC^2-C\\<\IM(\pi(w))-BC^2-Ci+C+BC^2-C=\IM(\pi(w))-Ci,
\end{multline*}
again by $t$-logarithmicity.

We claim that $\IM(\pi(\sum_{i=k-j}^kw_i))<\IM(\pi(w))-C(k-j-1)$.  The
proof is by induction on $j$.  If $j=0$, then the claim just says
$\IM(\pi(w_k))<\IM(\pi(w))-C(k-1)$, which is weaker than what
we already know.  Otherwise, we have
\begin{multline*}
\IM\left(\pi\left(\sum_{i=k-j}^kw_i\right)\right)=\IM\left(\pi\left(\sum_{i=k-j+1}^kw_i\right)+\pi(w_{k-j})\right)\\<\IM(\pi(w))-C(k-j)+C=\IM(\pi(w))-C(k-j-1),
\end{multline*}
where the inequality is by $t$-logarithmicity and induction. The claim
is proven.

Letting $j=k-1$ in the claim gives
\[
\IM(\pi(\sum_{i=k-j}^kw_i))=\IM(\pi(w))<\IM(\pi(w)),
\]
a contradiction.  The lemma follows if we let $D=BC^2$. Note that this
depends only on $A$ (via $C$) and $B$, not on $w$.
\end{proof}

We will want to extend the definition of $\phi$ to the new $F_S$, that
is the free group on the now arbitrary generating set $S$. The only
step which is not obvious is the definition of the map from $\sym{S}$
to $F_A\star\span{t}$. We simply choose the map once and for all,
requiring only that each element $s\in\sym{S}$ be mapped to a word
representing $s$ and that inverses be mapped to inverses. Then the
definition goes through without change.

If $w_1$, $w_2\in F_S$, then $\phi(w_1)$, $\phi(w_2)$ and
$\phi(w_1w_2)$ are all generalized $t$-words. If
$\func[\onto]{\alpha}{G}{\Z}$ is the projection, then
\[
\IM(\phi(w_1w_2))\le\max(\IM(\phi(w_1))+\alpha(w_2),\IM(\phi(w_2)))
\]
and
\[
\Im(\phi(w_1w_2))\ge\min(\Im(\phi(w_1))+\alpha(w_2),\Im(\phi(w_2))).
\]
Similarly, $\phi(w_1\inv)$ is a generalized $t$-word. We have
\[
\IM(\phi(w_1\inv))=\IM(\phi(w_1))-\alpha(w_1)
\]
 and
\[
\Im(\phi(w_1\inv))=\Im(\phi(w_1))-\alpha(w_1).
\]

Let $g\in G$. Let $i\in\Z$ and $k\in K$ be such that $g=t^ik$.  Define
$\IM(g)=\IM(k)$ and $\Im(g)=\Im(k)$.

\begin{cor}\label{fourth}
Let $K$ be a nontrivial abelian group and let $R$ be an automorphism
of $K$. Let $G=K\sd\span{t}$, where $t$ acts by $R$, so that $K$ has
the structure of a $\lp{t}$-module. Let $\func[\onto]{\alpha}{G}{\Z}$
be the projection. Suppose there is a finite $t$-efficient
$t$-logarithmic $t$-generating set for $K$. Let $S$ be any finite
generating set for $G$ and let $\func[\onto]{\sigma}{F_S}{G}$ be the
projection.  Let $z=\max\setst{\alpha(s)}{s\in\sym{S}}$.  For $w\in
F_S$, let $l(w)$ denote the length of $w$ as a word.  Then there is
$F\in\N$ such that, for every $w\in\ker(\alpha\circ\sigma)$, either
$\IM(\sigma(w))\le l(w)z/4+F$ or $\Im(\sigma(w))\ge-l(w)z/4-F$.
\end{cor}

\begin{proof}
Let $n=l(w)$. Let $w=s_1s_2\dots s_n$, $s_1$, \dots, $s_n\in\sym{S}$.
Consider the sequence of $n+1$ integers
\[
(0,\alpha(\sigma(s_n)),\alpha(\sigma(s_{n-1}s_n)),\dots,\alpha(\sigma(s_2\dots
s_n)),0=\alpha(\sigma(w))),
\]
where $\alpha$ again denotes the projection from $G$ to $\Z$. Either
at most $n/2$ of them are positive or at most $n/2$ are negative.
Assume without loss of generality that at most $n/2$ are positive. But
consecutive members of the sequence differ by at most $z$.  It follows
that, for each $i\in\W$, at most $n/2-2i$ are greater than $iz$.

Decompose $\phi(w)$ into the generalized $t$-words $v^+$ and $v^-$,
where $v^+$ consists of all letters of $\phi(w)$ coming from letters
$s_i$ of $w$ where
\[
\alpha(\sigma(s_{i+1}\dots s_n))>0
\]
and $v^-$ consists of all other letters of $\phi(w)$.  Let $N$ denote
the maximal length (as a generalized $t$-word) of the $\phi(s')$ for
all $s'\in\sym{S}$.  Let
\[
I=\max\setst{\IM(\phi(s'))}{s'\in\sym{S}}.
\]
Then the length of $v^-$ is at most $Nn$ and $\IM(v^-)\le I$.  Let $C$
be as in the definition of $t$-logarithmicity. Then $\IM(\pi(v^-))\le
I+C(\log_2(Nn)+1)$.

But, by the first paragraph and the definitions of $v^+$, $N$ and $I$,
for each $i\in\Z$, at most $N(n/2-2i)$ letters of $v^+$ are in
$\bigcup_{j=iz+I+1}^\infty t^jA$. Let $k=\gint{n/4-i+1}$ in the
preceding sentence. Then, for each $k\in\Z$, at most $2Nk$ letters of
$v^+$ are in $\bigcup_{j=nz/4+z+I+1-kz}^\infty t^jA$.

Apply Lemma~\ref{triangle} with $B=2N/z$.  Note that $B$ is
independent of $w$ and $n$. Then there are $D\in\N$ and $m\in\W$ also
independent of $w$ and $n$ such that more than $2Nm/z$ letters of
$v^+$ are in $\bigcup_{j=\IM(\pi(v^+))-D-m}^\infty t^jA$. Thus there
is $k=\gint{m/z}$ such that more than $2Nk$ letters of $v^+$ are in
\[
\bigcup_{j=\IM(\pi(v^+))-D-z-kz}^\infty t^jA.
\]
By the preceding paragraph, this is only possible if
\[
\IM(\pi(v^+))-D-z-kz<nz/4+z+I+1-kz,
\]
that is if $\IM(\pi(v^+))<nz/4+I+D+2z+1$.

But, by $t$-logarithmicity,
\begin{multline*}
\IM(\sigma(w))=\IM(\pi(\phi(w)))\\=\IM(\pi(v^+)+\pi(v^-))\le\max(\IM(\pi(v^+)),\IM(\pi(v^-)))+C\\\le\max\left(\frac{nz}{4}+I+D+2z+1,I+C(\log_2(Nn)+1)\right)+C<\frac{nz}{4}+F,
\end{multline*}
where $C$ is again as in the definition of $t$-logarithmicity and $F$
is a constant depending only, via $z$, $I$, $D$, $C$ and $N$, on $K$,
$R$ and $S$ (not on $n$ or $w$).
\end{proof}

We are now ready for the

\begin{proof}[Proof of Theorem~\ref{notac}]
Let $\alpha$ again denote the projection from $G$ to $\Z$. Let
$a\ne0\in A$ and let $s\in\sym{S}$ be chosen such that
$\alpha(s)\ge\alpha(s')$ for all $s'\in\sym{S}$.  Let $z=\alpha(s)$.
For ease of notation, let $C$ be the greater of the $C$s from the
definitions of $t$-efficiency and $t$-logarithmicity.  For $n\in\W$
and $i\in\Z$, let $g_n(i)=s^{n+i}as^{-2n}as^n=s^{i-n}as^{2n}as^{-n}$.
Then, for $\abs{i}\le n$, $\Abs[S]{g_n(i)}\le4n-\abs{i}+2\Abs[S]{a}$.

For $n\in\W$, define $h_n^+=g_n(J)$ and $h_n^-=g_n(-J)$, where
$J\in\W$ is a constant to be chosen later.  Let $n\ge J$. We have
\[
\Abs[S]{h_n^+}\le4n-J+2\Abs[S]{a}
\]
and
\[
\Abs[S]{h_n^-}\le4n-J+2\Abs[S]{a}.
\]
Since $h_n^+(h_n^-)\inv=s^{2J}$, we have
\[
\Abs[S]{h_n^+(h_n^-)\inv}\le2J.
\]
Note that this depends only on our choice of $J$. Let $N$ be the
length of $\phi(s)$ as a generalized $t$-word. Let
$v_n=\phi(s^{n+J}as^{-2n}as^n)=\phi(s^{n-J}as^{-2n}as^n)$. Then
$t^na+t^{-n}a$ is a $t$-word representing the same element as $v_n$,
and it is of minimal length so long as $n>C$.  We can ensure this by
picking $J>C$. Thus $\IM(v_n)-C=nz-C\le\IM(h_n^+)$.

Every edge of the (left) Cayley graph of $G$ with respect to $S$
connects some $g_1$ and $g_2\in G$ with
$\abs{\alpha(g_1)-\alpha(g_2)}\le z$.  Thus any path in the (left)
Cayley graph of $G$ with respect to $S$ connecting $h_n^-$ and $h_n^+$
must contain some $g\in G$ with $\abs{\alpha(g)}<z$.  Suppose
$\Abs[S]{g}\le4n-J+2\Abs[S]{a}$.  Then Corollary~\ref{fourth} says
that either $\IM(g)\le nz-Jz/4-\Abs[S]{a}z/2+F$ or $\Im(g)\ge
Jz/4+\Abs[S]{a}z/2-nz-F$, where $F$ is as in the corollary.  Without
loss of generality, we assume the former.

Since $A$ is $t$-logarithmic,
\[
\IM(h_n^+)\le\max(\IM(h_n^+g\inv)+\alpha(g),\IM(g))+C.
\]
We want to choose $J$ so that
\[
\IM(h_n^+)>\IM(g)+C.
\]
Since $\IM(g)+C\le nz-Jz/4-\Abs[S]{a}z/2+F+C$, it will suffice to take
\[
nz-\frac{Jz}{4}-\frac{\Abs[S]{a}z}{2}+F+C<nz-C,
\]
that is $J>(4/z)(F-\Abs[S]{a}z/2+2C)$. Note that this is independent
of $n$.  Then we will have $\IM(g)+C<nz-C\le\IM(h_n^+)$, as desired.
It will follow that
\[
nz\le\IM(h_n^+)+C\le\IM(h_n^+g\inv)+\alpha(g)+2C,
\]
that is $\IM(h_n^+g\inv)\ge nz-\alpha(g)-2C>nz-z-2C$.

Let $w\in F_S$ be of minimal length in $\sigma\inv(h_n^+g\inv)$.  For
$v$ a generalized $t$-word in $A$ and $i\in\Z$, let $v_i$ be the
number of letters of $v$ in $t^iA$. Let
$M=\sum_{i=-\infty}^\infty\max_{s'\in\sym{S}}\Abs[t^iA]{(\phi(s'))_i}$;
this sum is finite since $S$ is finite and, for each $s'\in\sym{S}$,
$\phi(s')$ has finitely many letters.  Let $B_{M,0}$ and $C_{M,0}$ be
as in the definition of strong $t$-logarithmicity.  As in the proof of
Corollary~\ref{fourth}, let $N$ denote the maximal length (that is
number of nonzero terms) of the $\phi(s')$ for all $s'\in\sym(S)$. Let
$I$ denote the maximal $\IM(\phi(s'))$ for all such $s'$. Then
\[
l(w)\ge\frac{\Abs[M]{\phi(w)}}{N}+2\frac{\IM(\phi(w))-I}{z}-\frac{\abs{\alpha(h_n^+g\inv)}}{z}.
\]

By strong $t$-logarithmicity and the result of two paragraphs ago,
\begin{multline*}
\Abs[M]{\phi(w)}>e^{\frac{\IM(h_n^+g\inv)-\IM(\phi(w))-C_{M,0}}{B_{M,0}}}-1\\>e^{\frac{nz-\IM(\phi(w))-z-2C-C_{M,0}}{B_{M,0}}}-1.
\end{multline*}
Also,
$\abs{\alpha(h_n^+g\inv)}\le\abs{\alpha(h_n^+)}+\abs{\alpha(g)}<J+z$.
Putting this all together, we have
\begin{multline*}
\Abs[S]{h_n^+g\inv}-2n\\>e^{\frac{nz-\IM(\phi(w))-z-2C-C_{M,0}}{B_{M,0}}}/N+\frac{2(\IM(\phi(w))-nz)}{z}-\frac{I+J}{z}-2.
\end{multline*}
This expression is bounded below as a function of $nz-\IM(\phi(w))$,
say by $F\in\Z$. Note that $F$ depends only on $A$, $S$ and the
one-time choices we made in defining the map $\phi$ (provided $J$ is
chosen appropriately). In particular, it is independent of $n$.  So
$\Abs[S]{h_n^+g\inv}\ge2n+F$.  Since $2n+F$ goes to infinity as $n$
does, we are done.
\end{proof}

\section{Bounded depth}\label{sbound}
This section is devoted to the proof of Theorem~\ref{bound}. To prove
this theorem, we will use a general lemma about groups obtained as the
semidirect product of an abelian group with $\Z$. We begin with the
following

\begin{defin}
Let $G$ be an indicable group and let $\func[\onto]{\phi}{G}{\Z}$.
Let $K=\ker\phi$. Fix a splitting $\alpha$ for $\phi$, so every
element of $G$ can be expressed uniquely as a product $\alpha(n)k$,
$k\in K$, $n\in\Z$. Let $A$ be a generating set for $G$. Then we call
$A$ \emph{symmetrized about $\Z$} if $\alpha(n)k\in A$ \tiff
$\alpha(n)k\inv\in A$ for all $k\in K$ and $n\in\Z$.
\end{defin}

\begin{notat}
If $A$ is symmetrized about $\Z$, then we have an involution on $A$,
which we denote with an overbar; thus
$\overline{\alpha(n)k}=\alpha(n)k\inv$.
\end{notat}

We extend the map $\overline{\cdot}$ to $A\inv$ so that it will be an
involution on $\sym{A}$.  We then extend it to $F_A$ so that
\[
\overline{a_1a_2\dots
a_m}=\overline{a_1}\overline{a_2}\dots\overline{a_m}
\]
for $m\in\N$ and $a_1$, \dots, $a_m\in\sym{A}$. 

\begin{lem}\label{invol}
Let $G$ be an indicable group and let $\func[\onto]{\phi}{G}{\Z}$ with
$K=\ker\phi$ abelian. Let $S$ be a generating set for $G$ symmetrized
about $\Z$ with respect to the splitting $\alpha$ and let
$\func[\onto]{\pi}{F_S}{G}$.  Let $w\in\pi\inv(K)$.  Then
$\pi(w)=\pi(\overline{w}\inv)$.
\end{lem}

\begin{proof}
Let $w=\alpha(n_1)k_1\dots\alpha(n_m)k_m$ with $n_1$, \dots,
$n_m\in\Z$, $k_1$, \dots, $k_m\in K$ and $\alpha(n_1)k_1$, \dots,
$\alpha(n_m)k_m\in\sym{S}$.  Then
$\pi(w)=\sum_{j=1}^mk_j^{\sum_{l=j+1}^mn_l}$, where we use additive
notation because $K$ is abelian. (We continue to denote the action of
$\Z$ by exponentiation in order to avoid confusion with the natural
action of $\Z$ on $K$ as an abelian group.)

On the other hand,
\[
\pi(\overline{w}\inv)=\pi(k_m\alpha(-n_m)\dots
k_1\alpha(-n_1))=\sum_{j=1}^mk_j^{-\sum_{l=1}^jn_l}=\sum_{j=1}^mk_j^{\sum_{l=j+1}^mn_l},
\]
where the last inequality is because $\sum_{l=1}^mn_l=0$ since
$w\in\pi\inv(K)$.  Since this is the same as $\pi(w)$, we are done.
\end{proof}

\begin{lem}\label{localityl}
Let $G$ be an indicable group and let $\func[\onto]{\phi}{G}{\Z}$ with
$K=\ker\phi$ abelian. Let $S$ be a generating set for $G$ symmetrized
about $\Z$ and let $\func[\onto]{\pi}{F_S}{G}$.  Then there exist $B$
and $C\in\N$ with the following property. Let $g\in G$ and $l\in\W$
with $l\le\phi(g)$.  Then there exist $w$, $w_1$ and $w_2\in F_S$ with
the following properties:
\begin{itemize}
\item $w=w_1w_2$ as words (that is without any cancellation),
\item $\pi(w)=g$,
\item $l(w)\le\Abs[S]{g}+C$,
\item $\abs{\phi(\pi(w_2))-l}$, $\abs{\phi(\pi(w_1))+l-\phi(g)}\le B$
  and
\item if $v$ is a prefix of $w_2$ or a suffix of $w_1$ then
  $\phi(\pi(v))\le2B$.
\end{itemize}
\end{lem}

\begin{proof}
Let $B=\max_{s\in S}\abs{\phi(s)}+1$.  Let $w'\in F_S$ with
$l(w')=\Abs[S]{g}$ and $\pi(w')=g$.  Let $T$ be the set of all $w''\in
F_S$ such that there exist $w_l$ and $w_r\in F_S$ with the following
properties:
\begin{itemize}
\item $w'=w_lw''w_r$ as words,
\item $\abs{\phi(\pi(w_r))-l}\le B$,
\item $\abs{\phi(\pi(w''w_r))-l}\le B$ and
\item there is no nonempty proper suffix $w_s$ of $w''$ with
  $\abs{\phi(\pi(w_sw_r))-l}\le B$.
\end{itemize}
Then $w'$ decomposes uniquely as $w_aw''_1\dots w''_kw_b$, where
\begin{itemize}
\item $k\in\W$,
\item $T=\setst{w''_i}{1\le i\le k}$,
\item $w_b$ is the minimal suffix of $w'$ with
  $\abs{\phi(\pi(w_b))-l}\le B$ and
\item $w_a$ is the minimal prefix with $\abs{\phi(\pi(w_a\inv
  w'))-l}\le B$.
\end{itemize}

Choose $j_0$, $j_1$, \dots, $j_k\in\pi\inv(\alpha(\set{-B,\dots,B}))$
such that
\[
\phi(\pi(j_k\inv w_b))=l
\]
and $\pi(j_{i-1}\inv w''_ij_i)\in K$.  Let $w'_a=w_aj_0$,
$w'''_i=j_{i-1}\inv w''_ij_i$ and $w'_b=j_k\inv w_b$.  Then
\begin{multline*}
\pi(w'_aw'''_1\dots w'''_kw'_b)=\pi(w_aj_0j_0\inv w''_1j_1j_1\inv\dots
w''_kj_kj_k\inv w_b)\\=\pi(w_aw''_1\dots w''_kw_b)=g,
\end{multline*}
and the $w'''_i$ commute with each other, since $K$ is abelian.

For $m$, $n\in\set{-B,\dots,B}$, let $w''_i\in T$ be an $(m,n)$-word
if $\pi(j_i)=\alpha(m)$ and $\pi(j_{i-1})=\alpha(n)$.  Let $W_{m,n}$
be the (possibly multi-) set of all $(m,n)$-words.  Let the finite
sequence $(w''_{m,n,i})$ be an enumeration of $W_{m,n}$.  Let
\[
w_{m,n}=j\inv\dots\overline{w''_{m,n,4}}\inv
w''_{m,n,3}\overline{w''_{m,n,2}}\inv w''_{m,n,1}j',
\]
where $j$, $j'\in\set{j_0,\dots,j_k}$ such that $\pi(j)=\alpha(n)$ and
$\pi(j')=\alpha(m)$.  Then
\[
\pi(w_{m,n})=\pi(\dots\overline{w'''_{m,n,4}}\inv
w'''_{m,n,3}\overline{w'''_{m,n,2}}\inv w'''_{m,n,1})=\pi(\dots
w'''_{m,n,2}w'''_{m,n,1}),
\]
by Lemma~\ref{invol}.

Let
\begin{equation}\label{newprod}
w=\overbrace{w'_a\prod_{0<m,n\le B}w_{m,n}\prod_{\abs{m},\abs{n}\le
    B,mn\le0}w_{m,n}}^{w_1}\overbrace{\prod_{-B\le
    m,n<0}w_{m,n}w'_b}^{w_2}.
\end{equation}
Then $\pi(w)=g$.

Let $J$ be at least the maximal length of any of the $j_i$; $J$
clearly can be taken independent of $g$. The length of $w_{m,n}$ is at
most $2J+$ the sum of the lengths of all the $(m,n)$-words $\in
W_{m,n}$. Then
\[
l(w)\le
l(w')+2\left[(2B+1)^2+1\right]J=\Abs[S]{g}+2\left[(2B+1)^2+1\right]J.
\]
Then we are done if we set $C=2\left[(2B+1)^2+1\right]J$ and let $w_1$
and $w_2$ be as marked in \eqref{newprod}. Note that $C$ is also
independent of $g$.
\end{proof}

\begin{cor}\label{locality}
Let $G$ be a group such that there is $\func[\onto]{\phi}{G}{\Z}$ with
$\ker\phi$ abelian. Let $S$ be a generating set for $G$ symmetrized
about $\Z$ and let $\func[\onto]{\pi}{F_S}{G}$.  Then there exist $B$
and $C\in\N$ with the following property. Let $g\in G$ with
$\phi(g)\ge0$.  Let $l_1$, $l_2\in\Z$ with $0\le l_1\le
l_2\le\phi(g)$.  Then there exists $w$, $w_1$, $w_2$ and $w_3\in F_S$
with the following properties:
\begin{itemize}
\item $w=w_1w_2w_3$ as words,
\item $\pi(w)=g$,
\item $l(w)\le\Abs[S]{g}+C$,
\item $\abs{\phi(\pi(w_3))-l_1}$, $\abs{\phi(\pi(w_2))+l_1-l_2}$ and
  $\abs{\phi(\pi(w_1))+l_2-\phi(g)}\le B$ and
\item if $v$ is a prefix of $w_2$ or $w_3$ or a suffix of $w_1$ or
  $w_2$ then $\phi(\pi(v))\le2B$.
\end{itemize}
\end{cor}

\begin{proof}
Let $w'$, $w'_1$ and $w'_2$ be the words $w$, $w_1$ and $w_2$ from
Lemma~\ref{localityl} with $l=l_1$. Let $w_r$ be the minimal suffix of
$w'$ such that $\phi(\pi(aw_r))>l_1+B$, where $a$ is the letter of
$w'$ next to the left of $w_r$.  Note that $w_r$ contains $w'_2$ by
construction, so we can let $w_r=w'_rw'_2$.  Let $w_l$, $w_{l1}$ and
$w_{l2}$ be the words $w$, $w_1$ and $w_2$ from Lemma~\ref{localityl}
with $g=\pi(w'w_r\inv)$ and $l=l_2-\phi(\pi(w_r))$.  We are done if we
let $B$ be as in Lemma~\ref{localityl}, $C$ be twice the $C$ from
Lemma~\ref{localityl}, $w=w_lw_r$, $w_1=w_{l1}$, $w_2=w_{l2}w'_r$ and
$w_3=w'_2$.
\end{proof}

For $g\in G$, as in Section~\ref{snotac}, we define
\[
\IM(g)=\IM(\alpha(\phi(g))\inv g)
\]
and
\[
\Im(g)=\Im(\alpha(\phi(g))\inv g).
\]
This makes sense since $\alpha(\phi(g))\inv g\in K$. As there, if $g$,
$h\in G$ and $C$ is as in the definition of $t$-logarithmicity,
\[
\IM(gh)\le\max(\IM(g)+\phi(h),\IM(h))+C,
\]
\[
\Im(gh)\ge\min(\Im(g)+\phi(h),\Im(h))-C,
\]
\[
\IM(g\inv)=\IM(g)-\phi(g)
\]
and
\[
\Im(g\inv)=\Im(g)-\phi(g).
\]
Furthermore, $\IM(g^h)=\IM(g)+\phi(h)$ and $\Im(g^h)=\Im(g)-\phi(h)$.

\begin{lem}\label{twoyes}
Let $G$ be an indicable group and $\func[\onto]{\phi}{G}{\Z}$ with
$K=\ker\phi$ abelian. Let $\func{\alpha}{\Z}{G}$ be the splitting map.
Let $t=\alpha(1)$, so that $t$ acts on $K$ by conjugation and makes it
into a $\lp{t}$-module. Let $A$ be a finite symmetrized $t$-efficient
$t$-logarithmic $t$-generating set for $K$, with $0\in A$. Let
\[
S=\setst{ta_1t^2a_2t}{a_1,a_2\in A}\cup A\cup\set{t}.
\]
For every $B\in\N$ there is $L\in\N$ such that the following holds.
Let $n\in\N$, $g$, $h$, $h_1$, $h_2$ and $h_3\in G$ such that
\begin{itemize}
\item $\phi(hg\inv)=4n$,
\item $\Im(hg\inv)\ge0$,
\item $\IM(hg\inv)\le4n$,
\item $h=h_1h_2h_3$,
\item $\abs{\phi(g\inv h_3)}\le B$,
\item $\abs{\phi(h_1)}\le B$,
\item $\Im(h_1)$, $\Im(h_2)\ge-B$,
\item $\IM(h_2)\le4n+B$ and
\item $\IM(h_3)\le\phi(g)+B$.
\end{itemize}
Then $\Abs[S]{g}\le\Abs[S]{h_1}+\Abs[S]{h_3}+2n+L$.
\end{lem}

\begin{proof}
Let $C$ be the constant with respect to which $A$ is
$t$-logarithmic. Let $w_1$, $w_2\in F_S$ be as given by
Lemma~\ref{localityl} with $l=\phi(g)$. Let $g_1=\pi(w_1)$ and
$g_2=\pi(w_2)$.  Then (possibly after increasing $B$)
\begin{itemize}
\item $g=g_1g_2$,
\item $\abs{\phi(g_1)}\le B$,
\item $\Im(g_1)\ge-B$ and
\item $\IM(g_2)\le\phi(g)+B$.
\end{itemize}

I claim that $\abs{\Abs[S]{g_2}-\Abs[S]{h_3}}\le F$. It follows from
$t$-logarithmicity and the conditions on the $h_i$ and $g_i$ that
\begin{multline*}
\Im(h_1h_2g_1\inv)\\\ge\min(\min(\Im(h_1)+\phi(h_2),\Im(h_2))-C-\phi(g_1),\Im(g_1)-\phi(g_1))-C\\\ge\min(\min(4n-3B,-B)-C-B,-2B)-C\\\ge\min(-4B-C,-2B)-C=-D,
\end{multline*}
where $D\in\N$ depends only on $A$ and $B$. But
\[
\abs{\phi(h_3g_2\inv)}\le\abs{\phi(g\inv h_3)}+\abs{\phi(g_1)}\le2B.
\]
Also,
\[
h_3g_2\inv=[(h_1h_2g_1\inv)\inv hg\inv]^{g_1}
\]
and $\phi(h_1h_2g_1\inv)=\phi(h)-\phi(h_3)-\phi(g_1)\le4n+2B$.  Thus,
we have
\begin{multline*}
\Im(h_3g_2\inv)=\Im((h_1h_2g_1\inv)\inv
hg\inv)+\phi(g_1)\\\ge\min(\Im(hg\inv),\Im(h_1h_2g_1\inv)-\phi(h_1h_2g_1\inv)+\phi(hg\inv))-C+\phi(g_1)\\\ge\min(0,-D-2B)-C-B\ge-E,
\end{multline*}
where $E\in\N$ also depends only on $A$ and $B$.  But, by
$t$-logarithmicity,,
\begin{multline*}
\IM(h_3g_2\inv)\le\max(\IM(h_3)-\phi(g_2),\IM(g_2)-\phi(g_2))+C\\\le\max(2B,2B)+C\le2B+C\le F,
\end{multline*}
where $F\in\N$ also depends only on $A$ and $B$.  Thus there is
$H\in\N$ (depending, via $C$, $D$, $E$, $F$ and the $t$-efficiency
constant, only on $A$ and $B$) such that $\Abs[S]{h_3g_2\inv}\le F$.
The claim follows.

For $j\in\set{1,2}$ and $i\in\Z$, let the $a_{ji}\in A$ be such that
\[
\alpha(\phi(g_1))\inv g_1=\sum_{i=-B}^\infty t^ia_{1(i+\phi(g_2))},
\]
\[
\alpha(\phi(g_2))\inv g_2=\sum_{i=-\infty}^{\phi(g)+B}t^ia_{2i}
\]
and these are minimal-length $t$-words representing
\[
\alpha(\phi(g_1))\inv g_1
\]
and
\[
\alpha(\phi(g_2))\inv g_2.
\]
(For ease of notation, let all $a_{ji}$ not referenced in the above
sums be $0$.) Then
\[
\sum_{i=1}^{4n-1}t^ia_{1(i+\phi(g_2))}
\]
and
\[
\sum_{i=1}^{4n-1}t^ia_{2i}
\]
are minimal-length $t$-words, so so are
\[
\sum_{i=1}^{4n-1}t^ia_{1(i+\phi(g))}
\]
and
\[
\sum_{i=1}^{4n-1}t^ia_{2(i+\phi(g))}.
\]
Let $J$ be the constant referred to as $C$ in the definition of
$t$-efficiency.  For each $n\in\N$, since $A$ is $t$-efficient and
$t$-logarithmic, there are $w_r(n)$ and $w_l(n)\in F_S$ such that
$l(w_r(n))$, $l(w_l(n))\le n+C+J+3$, $\phi(\pi(w_r(n)))=4n$,
$\phi(\pi(w_l(n)))=-4n$ and
\[
\pi(w_l(n)^tw_r(n))=\sum_{i=1}^{4n-1}t^i(a_{1(i+\phi(g))}+a_{2(i+\phi(g))}).
\]

Let $v_1\in\pi\inv(t^{-\phi(h_1)})$. For $i\in\Z$, $i\ge-B$, let
$h_{1i}\in A$ be such that $h_1=\sum_{i=-B}^\infty t^ih_{1i}$ and this
is a minimal-length $t$-word representing $h_1$; this is possible by
the assumption (in the hypothesis of the lemma) that
$\Im(h_1)\ge-B$. Let
\begin{multline*}
v_2\in\pi\inv\left(\sum_{i=0}^\infty
t^i(a_{1(i+\phi(g)+4n)}+a_{2(i+\phi(g)+4n)})\right.\\\left.{}-\sum_{i=-\phi(h_1)-B}^\infty
t^ih_{1(i+\phi(h_1))}\right).
\end{multline*}
By $t$-logarithmicity,
$\Im(\pi(v_2))\ge\min(0,-\phi(h_1)-B)-C\ge-2B-C$.  Although the
indices of both sums individually go to infinity, note that
\begin{multline*}
\pi(t^{\phi(g)+4n}v_2)=\sum_{i=\phi(g)+4n}^\infty
t^i(a_{1i}+a_{2i})-\sum_{i=-B}^\infty
t^{\phi(g)+4n-\phi(h_1)}h_{1i}\\=\alpha(\phi(g))\inv
g-\sum_{i=-\infty}^{\phi(g)+4n-1}t^i(a_{1i}+a_{2i})\\-t^{\phi(g)+4n-\phi(h_1)}\alpha(\phi(h_1))\inv
h_1+\sum_{i=\infty}^{1-B}t^{\phi(g)+4n-\phi(h_1)}h_{1i}\\=\alpha(\phi(g))\inv
g-t^{\phi(h_2h_3)}\alpha(\phi(h_1))\inv
h_1+P\\=t^{\phi(g)}\alpha(\phi(hg\inv))\inv
hg\inv+\alpha(\phi(h_2h_3))\inv h_2h_3+P,
\end{multline*}
where $P\in K$. Note that, by $t$-logarithmicity (and recalling that
$\abs{\phi(h_1)}\le B$), $\IM(P)\le\phi(g)+4n+1+2C$. It follows (again
by $t$-logarithmicity) that
\begin{multline*}
\IM(\pi(v_2))\\\le\max(\IM(hg\inv+\phi(g)),\IM(h_2h_3),\phi(g)+4n+1+2C)+2C-\phi(g)-4n\\\le\max(\phi(g)+4n+1+2C,\phi(g)+4n+2B+C)-\phi(g)-4n+2C\\\le\max(4C+1,2B+3C);
\end{multline*}
in particular, it has a bound depending only on $A$ and $B$ (not $n$).
Finally, let
\[
v_3\in\pi\inv\left(\alpha(\phi(g_1))\left[\sum_{i=-B}^{\phi(g_1)}t^ia_{1(i+\phi(g_2))}-\sum_{i=\phi(g_1)+1}^Bt^ia_{2(i+\phi(g_2))}\right]\right).
\]
Recall that $\abs{\phi(g_1)}\le B$ also.

Clearly,
\begin{multline*}
\phi(\pi(w_l(n)^t)\alpha(\phi(h_1))\inv
h_1^{\pi(v_1)}\pi(v_2w_r(n)v_3))\\=\phi(\pi(w_l(n)))-\phi(h_1)+\phi(h_1)+\phi(\pi(v_2))+\phi(\pi(w_r(n)))+\phi(\pi(v_3))\\=4n-\phi(h_1)+\phi(h_1)-4n+\phi(g_1)=\phi(g_1).
\end{multline*}
Also,
\begin{multline*}
\alpha(\phi(g_1))\inv\pi(w_l(n)^t)\alpha(\phi(h_1))\inv
h_1^{\pi(v_1)}\pi(v_2w_r(n)v_3)\\=\sum_{i=-B}^{\phi(g_1)}t^ia_{1(i+\phi(g_2))}-\sum_{i=\phi(g_1)+1}^Bt^ia_{2(i+\phi(g_2))}\\+\sum_{i=1}^{4n-1}t^{i+\phi(g_1)}(a_{1(i+\phi(g))}+a_{2(i+\phi(g))})\\+\sum_{i=0}^\infty
t^{i+\phi(g_1)+4n}(a_{1(i+\phi(g)+4n)}+a_{2(i+\phi(g)+4n)})\\-\sum_{i=-\phi(h_1)-B}^\infty
t^{i+\phi(g_1)+4n}h_{1(i+\phi(h_1))}+\sum_{i=-B}^\infty
t^{i+\phi(g_1)+4n-\phi(h_1)}h_{1i}\\=\sum_{i=-B}^{\phi(g_1)}t^ia_{1(i+\phi(g_2))}-\sum_{i=\phi(g_1)+1}^Bt^ia_{2(i+\phi(g_2))}\\+\sum_{i=\phi(g_1)+1}^{\phi(g_1)+4n-1}t^i(a_{1(i+\phi(g_2))}+a_{2(i+\phi(g_2))})\\+\sum_{i=\phi(g_1)+4n}^\infty
t^i(a_{1(i+\phi(g_2))}+a_{2(i+\phi(g_2))})\\-\sum_{i=\phi(g_1)+4n-\phi(h_1)-B}^\infty
t^ih_{1(i+\phi(h_1)-\phi(g_1)-4n)}\\+\sum_{i=\phi(g_1)+4n-\phi(h_1)-B}^\infty
t^ih_{1(i+\phi(h_1)-\phi(g_1)-4n)}\\=\sum_{i=-B}^\infty
t^ia_{1(i+\phi(g_2))}=\alpha(\phi(g_1))\inv g_1.
\end{multline*}
Thus
\[
\pi(w_l(n)^t)\alpha(\phi(h_1))\inv
h_1^{\pi(v_1)}\pi(v_2w_r(n)v_3)=g_1.
\]

For all $i\in\set{1,2,3}$,
\[
\abs{\phi(\pi(v_i))}\le B,
\]
\[
\IM(\pi(v_i))\le B+2C
\]
and
\[
\Im(\pi(v_i))\ge-2B-C.
\]
Thus there is $I\in\N$ (again depending only on $A$ and $B$) such that
all the $v_i$ are of length (as words in $F_S$) $\le I$.  Possibly
increasing $I$, we may arrange that $\Abs[S]{\alpha(\phi(h_1))}\le I$
as well; since $\abs{\phi(h_1)}\le B$ by hypothesis, $I$ still depends
only on $A$ and $B$. It follows that
\[
\Abs[S]{g_1}\le\Abs[S]{h_1}+5I+8+2n+2C+2J.
\]

Putting this all together, we get
\begin{multline*}
\Abs[S]{g}\le\Abs[S]{g_1}+\Abs[S]{g_2}\\\le\Abs[S]{h_1}+5I+8+2n+2C+2J+\Abs[S]{h_3}+F\le\Abs[S]{h_1}+\Abs[S]{h_3}+2n+L,
\end{multline*}
where we take $L=5I+8+2C+2J+F$.
\end{proof}

\begin{lem}\label{coef}
Let $G$ be an indicable group and $\func[\onto]{\phi}{G}{\Z}$ with
$K=\ker\phi$ abelian. Let $\alpha$ be the splitting map. Let $t$
represent the image of $1\in\Z$ under $\alpha$, so that $t$ acts on
$K$ by conjugation and makes it into a $\lp{t}$-module. Let $A$ be a
finite $t$-efficient $t$-logarithmic $t$-generating set for $K$. Let
$S$ be a finite generating set for $G$, symmetrized about $\Z$.  Let
$\func[\onto]{\pi}{F_S}{G}$.  Then there is $C\in\N$ with the
following property. Let $g\in G$ and let $w\in F_S$ be of minimal
length in $\pi\inv(g)$.  Let $v=\beta(w)\in N_A\ab$, where $\beta$ is
the same as $\phi$ in Section~\ref{snotac}.  Then for every $i\in\Z$
the length (with respect to $A^{t^i}$) of the component of $v$ in
$\Z^{A^{t^i}}$ is $\le C$.
\end{lem}

\begin{rk}
It follows that, under the above conditions, there is a bound on
$\abs{\IM(v)-\IM(g)}$ and $\abs{\Im(v)-\Im(g)}$.
\end{rk}

\begin{proof}
Since $S$ is finite, it will suffice to bound $n\in\N$ such that there
exist $s\in\sym{S}$ and $w_0$, $w_1$, \dots, $w_n\in F_S$ with
$w=w_0sw_1s\dots sw_{n-1}sw_n$ and
$\phi(\pi(w_1))=\dots=\phi(\pi(w_{n-1}))=-\phi(s)$.

Suppose $n$, $s$, $w_0$, \dots and $w_n$ are as above. We may assume
without loss of generality that $n$ is even; this will simplify the
notation.  We will bound $n$.  By Lemma~\ref{invol},
\[
w_0s\overline{s}\inv\overline{w_1}\inv
w_2s\overline{s}\inv\overline{w_3}\inv\dots
s\overline{s}\inv\overline{w_{n-1}}\inv w_n\in\pi\inv(g)
\]
also, and it is of minimal length since $w$ is. But
$s\overline{s}\inv\in K$ and, for any $i$ with $1\le i\le n/2-1$,
$\pi(\overline{w_{2i-1}}\inv w_{2i})\in K$, since
$\phi(\pi(w_{2i-1}))=\phi(\pi(w_{2i}))$.  Since $K$ is abelian,
\[
w_0\overline{w_1}\inv w_2\overline{w_3}\inv
w_4\dots\overline{w_{n-3}}\inv
w_{n-2}\left(s\overline{s}\inv\right)^{\frac{n}{2}}\overline{w_{n-1}}\inv
w_n\in\pi\inv(g)
\]
again, and again it is of minimal length since its length is the same
as that of $w$. But it follows from the $t$-logarithmicity of $A$
that, for sufficiently large $n$, $(s\overline{s}\inv)^{n/2}$ is not
of minimal length. The lemma is thus proven.
\end{proof}

We next prove the following proposition, from which
Theorem~\ref{bound} follows trivially.

\begin{prop}
Let $G$ be an indicable group and $\func[\onto]{\phi}{G}{\Z}$ with
$K=\ker\phi$ abelian. Let $\func{\alpha}{\Z}{G}$ be the splitting map.
Let $t=\alpha(1)$, so that $t$ acts on $K$ by conjugation and makes it
into a $\lp{t}$-module. Let $A$ be a finite symmetrized $t$-efficient
$t$-logarithmic $t$-generating set for $K$, with $0\in A$. Let
\[
S=\setst{ta_1t^2a_2t}{a_1,a_2\in A}\cup A\cup\set{t}.
\]
For all $i\in\N$, let $k_i$ be as in the statement of
Theorem~\ref{bound}.  Then there is $n\in\N$ with the following
property.  Let $g\in G$ with $\phi(g)\ge0$.  Let $C$ be the greatest
of
\begin{itemize}
\item $B$ from Lemma~\ref{localityl},
\item $B$ from Corollary~\ref{locality},
\item the $t$-logarithmicity constant,
\item the $t$-efficiency constant and
\item the bound from the remark following the statement of
  Lemma~\ref{coef}.
\end{itemize}
Let $h_g(n)\in G$, $v$, $v_1$ and $v_2$ $t$-words and $a_i$, $b_i\in
A$ for $i\in\Z$ have the following properties:
\begin{itemize}
\item $\phi(h_g(n))=\phi(g)+4n$,
\item $v$ represents $\alpha(\phi(h_g(n)g\inv))\inv h_g(n)g\inv$,
\item $\Im(v)\ge0$,
\item $\Im(v)\le4n$,
\item $v_1+v_2$ represents
\[
\alpha(\phi(h_g(n)))\inv h_g(n)-t^{\phi(g)+C-1}k_{4n-2C+2},
\]
\item $\IM(v_1)<\phi(g)+C$ and
\item $\Im(v_2)>\phi(g)+4n-C$.
\end{itemize}
Then $\Abs[S]{h_g(n)}>\Abs[S]{g}$.
\end{prop}

\begin{proof}
Let $w$, $w_1$ and $w_2\in F_S$ be the words given by
Lemma~\ref{localityl} with $l=\phi(g)$. For $n\in\N$, let $w(n)$,
$w_1(n)$, $w_2(n)$ and $w_3(n)\in F_S$ be the words given by
Corollary~\ref{locality} with $g=h_g(n)$, $l_1=\phi(g)$ and
$l_2=\phi(g)+4n=\phi(h_g(n))$.  We remind the reader of three salient
facts about these words, namely
\begin{itemize}
\item no prefix or suffix $w'$ of any $w_i$ or $w_i(n)$, except
  possibly a suffix of $w_2$ or $w_3(n)$ or a prefix of $w_1$ of
  $w_1(n)$, has $\phi(\pi(w'))<-2C$,
\item $\abs{\phi(\pi(w_3))-\phi(g)}$,
  $\abs{\phi(\pi(w_3(n)))-\phi(g)}\le C$ and
\item $\abs{\phi(\pi(w_2(n)))-4n}\le2C$.
\end{itemize}

For $i\in\set{0,1,2}$, let $m_i(n)$ be the number of letters of
$w_2(n)$ with $i$ of $a_1$, $a_2$ nontrivial, where $a_1$ and $a_2$
are as in the definition of $S$. Let $\beta$ be as in
Lemma~\ref{coef}.  Let $z(n)\in\N$ be the number of
$i\in\set{\phi(g),\dots,\phi(g)+4n-1}$ such that the coefficient of
$t^{i-\phi(w_3(n))}$ in $\beta(w_2(n))$ is $0$. Then
\[
4n-z(n)\le m_1(n)+2m_2(n)=2\Abs[S]{\pi(w_2(n))}-2m_0(n)-m_1(n)+2C,
\]
where for the latter inequality we remind the reader that $w_2(n)$ is
within $C$ of minimal length.  It follows that
$m_0(n)+m_1(n)\le2m_0(n)+m_1(n)\le2\Abs[S]{\pi(w_2(n))}-4n+z(n)+2C$.

The word $w_2(n)$ corresponds to a path between points of $\Z$
connecting $\phi(\pi(w_3(n)))$ with $\phi(\pi(w_2(n)w_3(n)))$. For
$i\in\set{0,1,2}$, we let $s_i(n)$ be the set of edges of this path
corresponding to letters $b\in\sym{S}$ with $i$ of $a_1$, $a_2$
nontrivial.  Thus $s_i(n)$ has cardinality $m_i(n)$. By a
\emph{stretch} we mean a maximal (under inclusion) set $T$ of adjacent
elements of
\[
\set{l+1,\dots,l+4n-1}
\]
such that
\begin{itemize}
\item $T$ lies entirely between $l+C$ and $l+4n-C$,
\item for every $m\in T$, every edge of the path incident to or
  passing over $m$ is $\in s_2(n)$ and
\item for every $m\in T$, the coefficient of $t^{m-\phi(\pi(w_3(n))}$
  in $\beta(w_2(n))$ is nontrivial.
\end{itemize}
The path must traverse each stretch an odd number of times, since it
begins at $\phi(\pi(w_3(n)))$ and ends at $\phi(\pi(w_2(n)w_3(n)))$.
But, if a stretch of $k$ integers is traversed at least three times,
this will take at least
\[
3\left(\frac{k}{4}-2\right)
\]
letters of $S$. If a stretch is traversed only once, then by the
second condition in the definition of a stretch the edges of this
traverse are only incident to integers of one parity, say even. Thus,
by the third condition in the definition of a stretch, the path must
also, if the stretch has any even integers in it, either enter the
stretch from one end, reach the integer one beyond the even integer in
it furthest from that end, and return, or else enter the stretch from
both ends and reach the same odd integer.  This will take at least
$2(k/4-2)$ letters, which, added to the $k/4-2$ letters consumed by
the one traverse, again makes $3(k/4-2)$ letters of $S$.

Let $N$ be the total number of integers in all stretches and $N_s$ the
number of stretches.  Then
\[
N\ge4n-z(n)-5(m_0(n)+m_1(n))-2C
\]
and
\[
N_s\le m_0(n)+m_1(n)+z(n)+1.
\]
Then we have (again using that $w_2(n)$ is within $C$ of having
minimal length)
\begin{multline*}
\Abs[S]{\pi(w_2(n))}\ge\frac{3N}{4}-6N_s-C\\\ge3n-\frac{7z(n)}{4}-\frac{39}{4}(m_0(n)+m_1(n))-\frac{5C}{2}-6\\\ge42n-\frac{39}{2}\Abs[S]{\pi(w_2(n))}-\frac{23z(n)}{2}-22C-6.
\end{multline*}
We thus get $E\in\N$ (depending, via $C$, only on $A$) such that
\[
\Abs[S]{\pi(w_2(n))}\ge\frac{84n}{41}-\frac{23z(n)}{41}-E.
\]

I claim there is $F\in\N$ independent of $n$ such that
\[
t^{\phi(g)+C-1}k_{4n-2C+2}-t^{\phi(\pi(w_3(n)))}\alpha(\phi(\pi(w_2(n))))\inv\pi(w_2(n))
\]
is represented by a generalized $t$-word of length at most $F$. To see
this, note that, by assumption, $v_1$ and $v_2$ are $t$-words such
that $-v_1-v_2$ represents
\begin{multline*}
t^{\phi(g)+C-1}k_{4n-2C+2}-t^{\phi(\pi(w_3(n)))}\alpha(\phi(\pi(w_2(n))))\inv\pi(w_2(n))\\-\alpha(\phi(\pi(w_3(n))))\inv\pi(w_3(n))\\-t^{\phi(\pi(w_2(n)w_3(n)))}\alpha(\phi(\pi(w_1(n))))\inv\pi(w_1(n)).
\end{multline*}
But, by assumption and Lemma~\ref{coef}, there is a $t$-word $v_3$
with $\IM(v_3)\le\phi(g)+2C$ representing
$\alpha(\phi(\pi(w_3(n))))\inv\pi(w_3(n))$.  Also, there is a $t$-word
$v_4$ with $\Im(v_4)\ge\phi(g)+4n-2C$ representing
\[
t^{\phi(\pi(w_2(n)w_3(n)))}\alpha(\phi(\pi(w_1(n))))\inv\pi(w_1(n)).
\]
Let $v_5$ and $v_6$ be minimal-length $t$-words representing the same
element as $v_3-v_1$ and $v_4-v_2$ respectively. Then
\[
\IM(v_5)\le\phi(g)+4C
\]
and
\[
\Im(v_6)\ge\phi(g)+4n-4C.
\]
Thus, for sufficiently large $n$, $v_5+v_6$ is a minimal-length (by
$t$-efficiency) $t$-word representing
\[
t^{\phi(g)+C-1}k_{4n-2C+2}-t^{\phi(\pi(w_3(n)))}\alpha(\phi(\pi(w_2(n))))\inv\pi(w_2(n)).
\]
But, by the first two conditions of Theorem~\ref{bound} and
Lemma~\ref{coef}, this means $\Im(v_5+v_6)\ge\phi(g)-4C$ and
$\IM(v_5+v_6)\le\phi(g)+4n+4C$. Since $C$ is independent of $n$, the
claim follows.

By the last condition of Theorem~\ref{bound}, it follows that there
are $I$, $J\in\N$ depending only on $A$ such that $23z(n)/41\le
I(\Abs[S]{\pi(w_2(n))}-2n)+J$.  Then
$\Abs[S]{\pi(w_2(n))}\ge84n/41+2In-I\Abs[S]{\pi(w_2(n))}-E-J$, so
\[
\Abs[S]{\pi(w_2(n))}\ge\frac{\frac{84n}{41}+2In-E-J}{1+I}\ge Kn-L,
\]
where $K>2$ and $L$ depend only on $A$.  Putting this together with
Lemma~\ref{twoyes} yields that
\begin{multline*}
\Abs[S]{g}\le\Abs[S]{\pi(w_1(n))}+\Abs[S]{\pi(w_3(n))}+2n+H\\\le\Abs[S]{h_g(n)}-\Abs[S]{\pi(w_2(n))}+2n+C+H\\\le\Abs[S]{h_g(n)}-(K-2)n+C+H+L,
\end{multline*}
where $H$ is the $L$ from the statement of Lemma~\ref{twoyes}, and
thus also depends only on $A$.  Rearranging yields
$\Abs[S]{h_g(n)}\ge\Abs[S]{g}+(K-2)n-C-H-L>\Abs[S]{g}$ for
sufficiently large $n$, where the definition of ``sufficiently large''
depends only on $A$.
\end{proof}

\section{Hyperbolic actions on abelian groups}\label{examples}
Suppose $K$ to be any nontrivial finite-rank torsion-free abelian
group, and let $t$ act by a hyperbolic automorphism $T$ whose matrix
has integer coefficients. We will construct a finite strongly
$t$-logarithmic $t$-generating set for $K$. Theorem~\ref{unbound} will
then give us a generating set for $K\sd\span{t}$ with respect to which
it has unbounded depth. Also, Theorem~\ref{notac} will tell us
$K\sd\span{t}$ is not almost convex for any generating set.

\subsection{Construction of the $t$-generating set}
To construct our $t$-generating set, we let $A$ be a basis for a
maximal-rank lattice in $K$. Then the $t$-generating set will be a
certain finite set of words in $A$. We will choose
this set of words so that it will be a $t$-generating set.

\begin{lem}\label{cover}
Let $K$ be a real finite-dimensional vector space and let $T$ be an
automorphism of $K$.  Suppose the (complex) eigenvalues of $K$ all
have absolute value $>1$.  Let $B$ be a basis of $K$ and let
$\overline{B}$ be the (closed) convex hull of $\sym{B}$.  Let
$\Abs{\cdot}$ be a norm on $K$.  Then there are $C_1>0$ and $C_2>1$
such that, for all $n\in\N$, $\overline{T^nB}$ contains the ball of
radius $C_1C_2^n$ with respect to $\Abs{\cdot}$.
\end{lem}

\begin{proof}
Since all norms on a finite-dimensional vector space are equivalent,
we may restrict attention to a norm such that $\Abs{T\inv}<1$; such a
norm exists since all (complex) eigenvalues of $T$ have absolute value
$>1$, so all eigenvalues of $T\inv$ have absolute value $<1$.

Let $C_1>0$ be such that $\overline{B}$ contains all $k\in K$ with
$\Abs{k}\le C_1$; this is possible since $B$ is a basis and $K$ is
finite-dimensional.  Let $n\in\N$ and $k\in K$ be such that
$\Abs{k}<C_1/\Abs{T\inv}^n$.  Then I claim that
$k\in\overline{T^nB}=T^n\overline{B}$.  This will be so \tiff
$T^{-n}k\in\overline{B}$.  But if $n\in\N$ we have
\[
\Abs{T^{-n}k}\le\Abs{T\inv}^n\Abs{k}<C_1.
\]
We are done by our choice of $C_1$ if we let $C_2=1/\Abs{T\inv}$; this
is $>1$ by our choice of norm $\Abs{\cdot}$.
\end{proof}

\begin{prop}\label{log}
Let $K$ be a finite-rank torsion-free abelian group and let $T$ be an
endomorphism of $K$. Then $T\ten\R$ is an endomorphism of $K\ten\R$.
Suppose it is an automorphism.  Suppose none of the (complex)
eigenvalues of $T$ have absolute value $1$.  Suppose further there is
a finite set $A\subseteq K$ such that
\begin{itemize}
\item $A$ is an $\R$-basis for $K\ten\R$.
\item $\span{TA}\subseteq\span{A}$ and
\item $\span{B}=K$, where $B=\bigcup_{i=-\infty}^\infty T^iA$.
\end{itemize}
Then (using additive notation in $K$) for every $k\in K$ there are
$C_1$ and $C_2\in\N$ such that, for all $n\in\N$, $\Abs[B]{nk}\le
C_1\log n+C_2$.
\end{prop}

\begin{proof}
Since none of the eigenvalues of $T$ have absolute value $1$,
$K\ten\R$ decomposes as the direct sum of an expanding subspace $K_e$
and an contracting subspace $K_c$. Let $A_e\subseteq K\ten\R$ denote
the projection of $A$ to $K_e$, and similarly let $A_c\subseteq
K\ten\R$ denote the projection of $A$ to $K_c$. For every $k\in K$,
let $k_e$ and $k_c\in K\ten\R$ denote the projections of $k$. Over
$\R$, $A_e$ spans $K_e$ and $A_c$ spans $K_c$. Let $d\in\N$ be the
maximal dimension of any generalized eigenspace. (If $K$ is trivial,
let $d=1$.)  Then trivially $B_e=\bigcup_{i=0}^{d-1}T^iA_e$ spans
$K_e$, and similarly $B_c=\bigcup_{i=0}^{d-1}T^iA_c$ spans $K_c$.

As in Lemma~\ref{cover}, we use the overbar to denote the symmetrized
closed convex hull of a set.  For $m\in\Z$, let
\[
E_m=\setst{k'\in
  K\ten\R}{k'_e\in\overline{T^mB_e},k'_c\in\overline{T^{-m}B_c}}.
\]
Let $\Abs{\cdot}$ denote the $\R$-norm with respect to $A$. It follows
from Lemma~\ref{cover} that there are $C_1>0$ and $C_2>1$ such that,
for every $k\in K$ and $m$, $n\in\N$, if $n\Abs{k}<C_1C_2^m$ then
$nk\in E_m$.  It follows that there are $D_1>0$ and $D_2>1$ such that,
for every $k\in K$ and $n\in\N$, there is $m\in\N$ such that $nk\in
E_m$ and $D_1D_2^m<n\Abs{k}$. We can choose $m$ for all $k$ and $n$ so
that $k\in\span{T^{-m}A}$; assume this done.

Let $D_e$ (respectively $D_c$) be the Hausdorff distance with respect
to $B_e$ (resp.\ $B_c$) between $\overline{TB_e}$
(resp.\ $\overline{T\inv B_c}$) and $\overline{B_e}$
(resp.\ $\overline{B_c}$). Let $D$ be the greater of $D_e$ and $D_c$.
Then any element of $\overline{T^mB_e}$ is within distance $D$ of
$\overline{T^{m-1}B_e}$ with respect to $T^{m-1}B_e$. Similarly, any
element of $\overline{T^{-m}B_c}$ is within distance $D$ of
$\overline{T^{1-m}B_c}$ with respect to $T^{1-m}B_c$. It follows that
any element $k\in E_m$ is within distance $D$ with respect to
$\bigcup_{i=0}^{d-1}T^{i+m-1}A$ of some $k'_1\in K$ such that
${k'_1}_e\in\overline{T^{m-1}B_e}$.  Similarly, $k$ is within distance
$D$ with respect to $\bigcup_{i=0}^{d-1}T^{i+1-m}A$ of some $k'_2\in
K$ such that ${k'_2}_c\in\overline{T^{1-m}B_c}$.

But all the eigenvalues of the restriction of $T$ to $K_c$ have
absolute value $<1$, and all the eigenvalues of its restriction to
$K_e$ have absolute value $>1$. Thus there is $F$ depending only on
$K$, $T$ and $A$ such that ${k'_1}_c$ is within $F$ of $k_c$ with
respect to $A_c$ and ${k'_2}_e$ is within $F$ of $k_e$ with respect to
$A_e$.  Let $k'=k'_1+k'_2-k$. Then $k'\in K$ is within $2D$ of $k$
with respect to $\bigcup_{i=0}^{d-1}T^{i+m-1}A\cup T^{i+1-m}A$.  Also,
$k'_e$ is within $F$ of $\overline{T^{m-1}B_e}$ with respect to $A_e$,
and $k'_c$ is within $F$ of $\overline{T^{1-m}B_c}$ with respect to
$A_c$.  We thus have $H$ depending only on $K$, $T$ and $A$ such that
$k'$ is within $H$ of $E_{m-1}$ with respect to $\Abs{\cdot}$. (Recall
that $\Abs{\cdot}$ is the $\R$-norm with respect to $A$.)

Repeating this process $m$ times, we find that any $k\in E_m$ is
within $2mD$ with respect to $\bigcup_{i=1-m}^{m+d-2}T^iA$ of some
$k''\in K$ within $mH$ of $E_0$ with respect to $\Abs{\cdot}$.  Let
$I$ denote the radius of $E_0$ with respect to $\Abs{\cdot}$. Then
$\Abs{k''}\le I+mH$.  Thus, by another application of
Lemma~\ref{cover}, there exist $D_3$, $D_4\in\N$ depending only on
$K$, $T$ and $A$ such that $k'\in E_{D_3\log m+D_4}$.  For $m$ large,
$D_3\log m+D_4<m$.  Thus, we can repeat the procedure in this and the
preceding two paragraphs to find $D_5$, $D_6\in\N$ depending only on
$K$, $T$ and $A$ such that, for all $m\in\N$, any $k\in E_m$ is, with
respect to $\bigcup_{i=1-m}^{m+d-2}T^iA\subseteq B$, within $D_5m$ of
some $k'''\in K$ with $\Abs{k'''}<D_6$.

Putting the preceding paragraph together with the second paragraph
tells us that there exist $D_1$, $D_2$, $D_5$, $D_6\in\N$ depending
only on $K$, $T$ and $A$ with the following property. Let $k\in K$ and
$n\in\N$.  Then there is $m\in\N$ such that $D_1D_2^m<n\Abs{k}$.
Also, $nk$ is, with respect to $\bigcup_{i=1-m}^{m+d-2}T^iA\subseteq
B$, within distance $D_5m$ of some $k'_n\in K$ with $\Abs{k'_n}<D_6$.
Finally, $k'_n\in\span{nk,T^{1-m}A}=\span{T^{-m}A}$, since
$\span{TA}\subseteq\span{A}$.

I claim there is $L\in\N$ with the following property. Let
$a\in\span{A}$ with $\Abs{a_c}\le D_6$. Then there is $a'\in\span{TA}$
within $H$ of $a$ with respect to $A$ and with $\Abs{a'_c}\le D_6$.
To prove this claim, let $a\in\span{A}$ with $\Abs{a_c}\le D_6$. If
$K_e$ is trivial, then $a_c=a$, so $\Abs{a}\le D_6$, and, for
sufficiently large $i$, we have $\Abs{T^ia}<\Abs{a}\le D_6$.  Since
$T^ia\in\span{TA}$ for $i\ge1$, we can assume $K_e$ is not trivial.
There is $J$ depending only on $K$, $T$ and $A$ such that, for every
$a$, there is $a'\in\span{TA}$ within $J$ of $a$ with respect to $A$.
Then $\Abs{a'_c}\le D_6+J$.  If $K_c$ is trivial, then in particular
$\Abs{a'_c}=0\le D_6$, which is what we want, so suppose $K_c$ is also
not trivial.  The projection of $\span{TA}$ is dense in $K_c$ since
$T\ten\R$ is an automorphism. Thus there is $a''\in\span{TA}$
depending only on $K$, $T$ and $A$ with $0<\Abs{a''_c}<2D_6$.  Then
there is $i\in\Z$ with $\abs{i}\le(J+D_6)/\Abs{a''_c}+1$ such that
$\Abs{a'_c+ia''_c}\le D_6$.  The claim is proven.

Applying the claim $m$ times to the situation of the preceding
paragraph, we see that there exist $D_1$, $D_2$, $D_5$, $D_6$,
$D_7\in\N$ depending only on $K$, $T$ and $A$ with the following
property.  Let $k\in K$ and $n\in\N$.  Then there is $m\in\N$ such
that $D_1D_2^m<n\Abs{k}$.  Also, $nk$ is, with respect to
$\bigcup_{i=1-m}^{m+d-2}T^iA\subseteq B$, within distance $D_5m$ of
some $k'_n\in K$ with $\Abs{k'_n}<D_6$. Finally, $k'_n$ is, with
respect to $\bigcup_{i=-m}^{-1}T^iA\subseteq B$, within distance
$D_7m$ of some $k''_n\in\span{A}$ such that $\Abs{{k''_n}_c}\le D_6$.

But, since $A$ is finite, there is $M\in\N$ such that, for all
$i\in\N$ and $a\in A$, $\Abs{(T^{-i}a)_e}<M$.  It follows that
$\Abs{(k'_n-k''_n)_e}\le D_7Mm$.  Thus there are $D_8$ and $D_9\in\N$
depending only on $K$, $T$ and $A$ such that
\[
\Abs{(k''_n)_e}\le\Abs{(k'_n)_e}+\Abs{(k'_n-k''_n)_e}<D_8m+D_9.
\]
Thus there is $D_{10}\in\N$ such that
$\Abs{k''_n}\le\Abs{(k''_n)_e}+\Abs{(k''_n)_c}<D_8m+D_{10}$.  Recall
that $\Abs{\cdot}$ denotes the norm with respect to $A\subseteq B$.
Thus, by the triangle inequality, there is $D_{11}\in\N$
($=D_5+D_7+D_8$) depending only on $K$, $T$ and $A$ such that $nk$ is
within $D_{11}m+D_{10}$ of $0$ with respect to $B$, so we are done.
\end{proof}

The following is an easy consequence.

\begin{cor}\label{lowerc}
Let $K$ be a finite-rank torsion-free abelian group and let $T$ be an
endomorphism of $K$. Then $T\ten\R$ is an endomorphism of $K\ten\R$.
Suppose it is an automorphism.  Suppose none of the (complex)
eigenvalues of $T$ have absolute value $1$.  Suppose finally there is
a finite set $A\subseteq K$ such that
\begin{itemize}
\item $A$ is an $\R$-basis for $K\ten\R$.
\item $\span{TA}\subseteq\span{A}$ and
\item $\span{B}=K$, where $B=\bigcup_{i=-\infty}^\infty T^iA$.
\end{itemize}
Let $\func[\onto]{\pi}{\Z^B}{K}$ be the projection. Then there is
$n\in\N$ such that, for any $k\in K$ and $w\in\Z^B$ of minimal length
in $\pi\inv(k)$, the coefficient of every letter of $B$ in $w$ has
absolute value less than $n$. In particular, $\setst{\sum_{a\in
    A}i_aa}{\abs{i_a}<n,a\in A}$ is a $t$-generating set for $K$, if
$t$ acts by $T$.
\end{cor}

\subsection{Proof of strong $t$-logarithmicity}
We will next show the $t$-generating set just constructed is strongly
$t$-logarithmic.  We will use that, roughly speaking, each element of
$K$ is represented by finitely many minimal-length $t$-words.  More
precisely, we have

\begin{prop}\label{finite}
Let $K$ be a finite-rank torsion-free abelian group and let $L$ be a
full-rank lattice in $K$. Let $T$ be a hyperbolic automorphism of $K$
which acts on $L$ by an endomorphism and gives $K$ the structure of a
$\lp{t}$-module.  Let $A$ be a finite $t$-generating set for $K$.
Then for every $k\in K-\set{0}$ and $l\in\N$ there are only finitely
many minimal $t$-words $w$ in $A$ of length at most $l$ representing
$k$.
\end{prop}

The proof is deferred to Subsection~\ref{ssfinite}.

\begin{prop}\label{linlog}
Let $K$ be a finite-rank torsion-free abelian group and let $T$ be a
hyperbolic automorphism of $K$, so that $K$ has the structure of a
$\lp{t}$-module.  Choose a norm on $K\ten\R$ and let $d_c$
(respectively $d_e$) denote distance with respect to this norm from
the contracting (resp.\ expanding) subspace. Let $A$ be a finite
$t$-generating set for $K$.  For $w$ a $t$-word in $A$, let $\pi(w)$
be the element of $K$ it represents. Let $L\supseteq A$ be a full-rank
lattice in $K$ on which $T$ acts by an endomorphism.  Then
$\sum_{i=1}^\infty\max\setst{d_e(T^ia)}{a\in A}$ is finite.  Also, for
every $D\in\R$ and $n\in\W$, there is $E_1\in\Z$ with the following
property.  Let $w\in\pi\inv(L)$ be a $t$-word.  Suppose
$d_c(\pi(w))\le D$ and $w$ is minimal and within $n$ of minimal length
in $\pi\inv(\pi(w))$ subject to the condition $\Im(w)\ge0$.  Then
$\IM(w)<E_1$.

Similarly, $\sum_{i=1}^\infty\max\setst{d_c(T^{-i}a)}{a\in A}$ are
finite.  Also, for every $D\in\R$ and $n\in\W$, there is $E_2\in\Z$
with the following property.  Let $w\in\pi\inv(L)$ be a $t$-word.
Suppose $d_e(\pi(w))\le D$ and $w$ is minimal and within $n$ of
minimal length in $\pi\inv(\pi(w))$ subject to the condition
$\IM(w)\le0$.  Then $\Im(w)>E_2$.
\end{prop}

\begin{proof}
If $w=0$, then any $E_1$ and $E_2$ at all will work. We thus may
assume $w\ne0$. It follows that $w\notin\ker\pi$, since the only
minimal element of $\ker\pi$ is $0$.

We prove the first paragraph of the proposition; the proof of the
second is exactly analogous.

Since $T$ is hyperbolic, there are $C_1>0$ and $C_2<1$ such that, for
every $k\in K$, $d_e(T^ik)<C_1C_2^id_e(k)$.  In particular, this is so
for every $k\in A$.  Thus
\[
D'=\sum_{i=1}^\infty\max\setst{d_e(T^ia)}{a\in A}
\]
is finite, proving the first claim.

For the second claim, let $w$ be a $t$-word as in the condition of the
proposition.  Since $\Im(w)\ge0$, it follows from the first claim that
$d_e(\pi(w))<D'$.  Then $d_e(\pi(w))<D'$ and $d_c(\pi(w))\le D$. Thus
$\pi(w)$ is constrained to lie within a bounded region.  Since
$\pi(W)\in L$ and $L$ is discrete, there are finitely many possible
values of $\pi(w)$. None of these is $0$, since we are assuming
$w\notin\ker\pi$.  By Lemma~\ref{finite}, each of them is represented
by finitely many minimal $t$-words $w$ of any given length, in
particular of any length within $n$ of that length which is least
subject to the condition that $\Im(w)\ge0$. Thus there are finitely
many possibilities for $w$. We are done if we let $E_1$ be the maximum
of their $\IM$'s.
\end{proof}

\begin{prop}\label{upper}
Let $K$ be a finite-rank torsion-free abelian group and let $T$ be a
hyperbolic automorphism of $K$, so that $K$ has the structure of a
$\lp{t}$-module.  Choose a norm on $K\ten\R$ and let $d_c$
(respectively $d_e$) denote distance with respect to this norm from
the contracting (resp.\ expanding) subspace.  Let $A$ be a
$t$-generating set for $K$.  Let $L\supseteq A$ be a full-rank lattice
in $K$ on which $T$ acts by an endomorphism. For $w$ a $t$-word, let
$\pi(w)$ denote the element of $K$ it represents. Then
$d_1=\sum_{i=1}^\infty\max\setst{d_c(T^{-i}a)}{a\in A}$ and
$d_2=\sum_{i=1}^\infty\max\setst{d_e(T^ia)}{a\in A}$ are finite.  Let
$N_1$ and $N_2$ be the radius-$d_1$ and -$d_2$ closed neighborhoods,
respectively, of the contracting (resp.\ expanding) subspace of
$K\ten\R$.  Let $n\in\W$. Then there are $E_1$ and $E_2\in\Z$ such
that, for any $k\in K\cap N_1$ (resp.\ $L\cap N_2)$, any minimal
$t$-word $w$ representing $k$ and within $n$ of minimal length in
$\pi\inv(\pi(w))$ has $\IM(w)<E_1$ (resp.\ $\Im(w)>E_2$).  In
particular, $\IM(k)<E_1$ (resp.\ $\Im(k)>E_2$), for $E_1$ and $E_2$
independent of $n$.
\end{prop}

\begin{proof}
The sums $d_1$ and $d_2$ are finite by Proposition~\ref{linlog}.

Let $E_1$ and $E_2$ be given by Proposition~\ref{linlog} with
$D=2\max(d_1,d_2)$ and $n$ as in this proposition.  Let $w$ be a
minimal $t$-word within $n$ of minimal length in $\pi\inv(\pi(w))$ and
such that $\IM(w)\ge E_1$ (resp.\ $w\in\pi\inv(L)$ and $\Im(w)\le
E_2$).  We will show that $\pi(w)\notin N_1$ (resp.\ $N_2$). Break $w$
up into two parts, $w_h$ and $w_t$, where $w_h$ consists of the terms
from $\bigcup_{i=0}^\infty T^iA$ (resp.\ $\bigcup_{i=-\infty}^0T^iA$)
and $w_t$ of those from $\bigcup_{i=-\infty}^{-1}T^iA$
(resp.\ $\bigcup_{i=1}^\infty T^iA$). Thus $w_h$ and $w_t$ are minimal
$t$-words, and $w=w_h+w_t$.  It follows that
$\pi(w)=\pi(w_h)+\pi(w_t)$.  Then $\pi(w_h)\in L$ (in the first case
trivially, in the second since $\pi(w)$ and $\pi(w_t)$ are both $\in
L$).  Furthermore, $\Im(w_h)\ge0$ (resp.\ $\IM(w_h)\le0$) and
$\IM(w_h)=\IM(w)\ge E_1$ (resp.\ $\Im(w_h)\le E_2$).

Since $w$ is only within $n$ of minimal length among $t$-words, $w_h$
is not necessarily within $n$ of minimal length in
$\pi\inv(\pi(w_h))$.  However, $w_h$ is within $n$ of minimal length
in $\pi\inv(\pi(w_h))$ subject to the condition that $\Im(w_h)\ge0$
(resp.\ $\IM(w_h)\le0$).  Thus, by Proposition~\ref{linlog},
$d_c(\pi(w_h))>D\ge2d_1$ (resp.\ $d_e(\pi(w_h))>D\ge2d_2$).  But
clearly $d_c(\pi(w)-\pi(w_h))=d_c(\pi(w_t))\le d_1$
(resp.\ $d_e(\pi(w)-\pi(w_h))\le d_2$) since $w_t$ is a $t$-word. By
the triangle inequality, it follows that $d_c(\pi(w))>d_1$
(resp.\ $d_e(\pi(w))>d_2$). Thus $k\notin N_1$ (resp.\ $N_2$), as
claimed.

The last sentence follows trivially.
\end{proof}

\begin{cor}\label{upperc}
Let $K$ be a finite-rank torsion-free abelian group and let $T$ be a
hyperbolic automorphism of $K$, so that $K$ has the structure of a
$\lp{t}$-module.  Choose a norm on $K\ten\R$ and let $d_c$
(respectively $d_e$) denote distance with respect to this norm from
the contracting (resp.\ expanding) subspace.  Let $A$ be a
$t$-generating set of $K$.  Let $L\supseteq A$ be a full-rank lattice
in $K$ on which $T$ acts by an endomorphism. For $w$ a $t$-word, let
$\pi(w)$ denote the element of $K$ it represents. Then
$d_1=\sum_{i=1}^\infty\max\setst{d_c(T^{-i}a)}{a\in A}$ and
$d_2=\sum_{i=1}^\infty\max\setst{d_e(T^ia)}{a\in A}$ are finite.  Let
$N_1$ and $N_2$ be the radius-$d_1$ and -$d_2$ closed neighborhoods,
respectively, of the contracting (resp.\ expanding) subspace of
$K\ten\R$.  Let $n\in\W$. Then there are $E_1$ and $E_2\in\Z$ such
that, for any $i\in\Z$ and any $k\in K\cap T^iN_1$ (resp.\ $T^iL\cap
T^iN_2)$, any minimal $t$-word $w$ representing $k$ and within $n$ of
minimal length in $\pi\inv(\pi(w))$ has $\IM(w)<E_1+i$
(resp.\ $\Im(w)>E_2+i$).  In particular, $\IM(k)<E_1+i$
(resp.\ $\Im(k)>E_2+i$), for $E_1$ and $E_2$ independent of $n$ (and
$i$).
\end{cor}

\begin{proof}
Just apply $T^i$ to the statement of Proposition~\ref{upper}.
\end{proof}

\begin{prop}\label{strlog}
Let $K$ be a nontrivial finite-rank torsion-free abelian group and $T$
a hyperbolic automorphism of $K$. Let $L$ be a full-rank lattice in
$K$ on which $T$ acts by an endomorphism. Suppose
$\bigcup_{i=-\infty}^\infty T^iL=K$. Let $G=K\sd\span{t}$, where $t$
acts by $T$, so that $K$ has the structure of a $\lp{t}$-module.  Then
$K$ has a finite strongly $t$-logarithmic $t$-generating set.
\end{prop}

\begin{proof}
Let $A$ be a basis of $L$. Then the hypotheses of
Corollary~\ref{lowerc} are clearly satisfied, so let $A'$ be the
finite $t$-generating set given by that result. We will show $A'$ is
strongly $t$-logarithmic.

Choose a norm on $K\ten\R$ and let $d_c$ (respectively $d_e$) denote
distance with respect to this norm from the contracting
(resp.\ expanding) subspace. Let $N_1$, $N_2$, $E_1$ and $E_2$ be
defined as in Corollary~\ref{upperc} with $A=A'$ and with respect to
the above norm.  Let $B=\bigsqcup_{i=-\infty}^\infty T^iA'$ and let
$\func[\onto]{\pi}{\Z^B}{K}$ be the projection.  Let $m\in\N$, $k\in
K$ and $w\in\pi\inv(k)-\set{0}$. (Note that $w$ is a generalized
$t$-word, in the terminology of Section~\ref{defs}, but not
necessarily a $t$-word.) Let $\Abs[m]{w}$ be as in the definition of
strong $t$-logarithmicity. Then
\begin{multline*}
d_c(T^{-\IM(w)}k)\\\le m\sum_{i=0}^\infty\max\setst{d_c(T^{-i}a)}{a\in
A'}+\Abs[m]{w}\max_{i=0}^\infty\max\setst{d_c(T^{-i}a)}{a\in
A'}\\\le(m+\Abs[m]{w})\sum_{i=0}^\infty\max\setst{d_c(T^{-i}a)}{a\in
A'},
\end{multline*}
where the first inequality is by the definition of $\Abs[m]{w}$. Thus
\[
k\in(m+\Abs[m]{w})T^{\IM(w)}N_1.
\]
It follows that there are $C_1>0$ and $C_2\in\R$ (depending only on
$K$ and $T$) such that
\[
k\in T^{\IM(w)+C_1\log(m+\Abs[m]{w})+C_2}N_1.
\]

Let $w'$ be a minimal $t$-word representing $k$ and within $n$ of
minimal length among all $t$-words representing $k$. (The set of all
such $t$-words is $\pi\inv(\pi(w))$ in the notation of
Corollary~\ref{upperc}.)  Let $E_1$ and $E_2$ be given by
Corollary~\ref{upperc} with $n$ as in this proof. Then, by
Corollary~\ref{upperc},
\[
\IM(w')<\IM(w)+C_1\log(m+\Abs[m]{w})+C_2+E_1.
\]

In exactly the same way,
\[
k\in T^{\Im(w)-C_1\log(m+\Abs[m]{w})+C_2}N_2.
\]
Also, since every letter of $w$ is $\in T^{\Im(w)}L$, $k\in
T^{\Im(w)}L$.  By Corollary~\ref{upperc} again,
\[
\Im(w')>\Im(w)-C_1\log(m+\Abs[m]{w})+C_2+E_2.
\]
\end{proof}

\begin{cor}
Let $K$ be a nontrivial finite-rank torsion-free abelian group and $T$
a hyperbolic automorphism of $K$.  Let $L$ be a full-rank lattice in
$K$ on which $T$ acts by an endomorphism.  Suppose
$\bigcup_{i=-\infty}^\infty T^iL=K$. Let $G=K\sd\span{t}$, where $t$
acts by $T$.  Then $G$ has a generating set with respect to which it
has deep pockets. Furthermore, it is not almost convex with respect to
any generating set.
\end{cor}

\begin{proof}
Combine Theorems~\ref{unbound}, \ref{notac} and
Proposition~\ref{strlog}.
\end{proof}

\subsection{Proof of Proposition~\ref{finite}}\label{ssfinite}
In this subsection, we complete the proof of Proposition~\ref{strlog}
by proving Proposition~\ref{finite}, whose statement we repeat here
for the convenience of the reader.

\begin{rest}{Proposition~\ref{finite}}
Let $K$ be a finite-rank torsion-free abelian group and let $L$ be a
full-rank lattice in $K$. Let $T$ be a hyperbolic automorphism of $K$
which acts on $L$ by an endomorphism and gives $K$ the structure of a
$\lp{t}$-module.  Let $A$ be a finite $t$-generating set for $K$.
Then for every $k\in K-\set{0}$ and $l\in\N$ there are only finitely
many minimal $t$-words $w$ in $A$ of length at most $l$ representing
$k$.
\end{rest}

The plan of the proof is that $K$ (or at least a finite-index
submodule of $K$) can be split as a direct sum of two pieces: a
finitely generated piece, which we call $K_d$, and a piece $K_c$ with
$\bigcap_{i=-\infty}^\infty T^ik_c=\set{0}$. This splitting is given
by Lemma~\ref{gsplit}. Then we will deal with $K/K_c$ in
Lemma~\ref{find} and with $K/K_d$ in Corollary~\ref{finccc}.

\begin{lem}[splitting]\label{gsplit}
Let $K$ be a finite-rank torsion-free abelian group and let $T$ be a
hyperbolic automorphism of $K$. Let $L$ be a full-rank lattice in $K$
on which $T$ acts as an endomorphism. Suppose
$\bigcup_{i=-\infty}^\infty T^iL=K$.  Let
$K_d=\bigcap_{i=-\infty}^\infty T^iL$. Then $K$ has a finite-index
$T$-submodule $K'$ such that $K_d$ is a complemented $T$-submodule of
$K'$.  Furthermore, $K_d$ has a complement $K_c$ in $K'$ such that
$K/K_d$ and $K/K_c$ are both torsion-free.
\end{lem}

\begin{proof}
It is clear that the actions by $T$ and $T\inv$ preserve $K_d$. It
remains to show that it is complemented in some finite-index submodule
of $K$.

First, I claim that $K/K_d$ is torsion-free. To this end, let $L_d$
denote the set of all $l\in L$ such that there is $n\in\N$ with $nl\in
K_d$.  Then $K_d\subseteq L_d$ and, for all $i\in\W$, $T^iL_d\subseteq
L_d$.  Let $i\in\W$, $l\in L$ and $T^il\in L_d$.  Then there is some
$n\in\N$ with $nT^il=T^i(nl)\in K_d$.  It follows that $nl\in K_d$, so
$l\in L_d$. Thus also $T^il\in L_d$.  We have thus shown that
$T^iL\cap L_d\subseteq T^iL_d$. Since clearly $T^iL_d\subseteq T^iL$,
we have $T^iL\cap L_d=T^iL_d$.  Thus
\[
K_d=K_d\cap L_d=\bigcap_{i=0}^\infty(T^iL\cap
L_d)=\bigcap_{i=0}^\infty T^iL_d.
\]
But $K_d$ must be a full-rank subgroup of $L_d$, by the definition of
$L_d$.  Since $L$ is finitely generated, so is $L_d$, so $L_d/K_d$ is
finite.  By the above equation, this implies there is $i\in\N$ such
that $T^iL_d=K_d$.  Now let $k\in K-K_d$ and $n\in\N$ such that $nk\in
K_d$.  Let $j\in\Z$ be such that $k\in T^jL$. Then $T^{i-j}k\in
T^iL-K_d$, but $nT^{i-j}k=T^{i-j}(nk)\in K_d$. Thus $T^{i-j}k\in
T^iL_d-K_d$, a contradiction. Our claim is proven.

Thus $K_d$ and $K/K_d$ are both finite-rank torsion-free abelian
groups.  In fact, $K_d$ is finitely generated since $L$ is, so the
action of $T$ on $K_d$ has all its (complex) eigenvalues algebraic
units.  Suppose some eigenvalue of $T$ on $K/K_d$ were a unit. Then
let $K_u$ be the subspace of $K/K_d\ten\Q$ generated by the
generalized eigenspaces of that eigenvalue and all its conjugates over
$\Q$.  But then $K_u\cap L/K_d$ is a nontrivial finitely generated
subgroup of $K/K_d$ on which $T$ acts as an automorphism. This is a
contradiction, so none of the eigenvalues of the action of $T$ on
$K/K_d$ are units.

It follows by rational canonical form that $K_d\ten\Q$ is complemented
as a $T$-submodule of $K\ten\Q$.  Let $K_c$ be the intersection of the
complement with $K$. Clearly, $K_c$ is a $T$-module, and all elements
of $K$ which have a multiple in $K_c$ are themselves in $K_c$.  Thus
$K/K_c$ is torsion-free.

Let $K'=\span{K_c,K_d}$; we will show that $\ind{K}{K'}<\infty$. Note
that $K'$ is of full rank in $K$, since $K\ten\Q=K_d\ten\Q\ds
K_c\ten\Q$.  Thus $\ind{L}{K'\cap L}<\infty$, since $L$ is finitely
generated.  Since $K'$ is invariant under $T$ and $T\inv$ (since $K_d$
and $K_c$ are), it follows that, for all $i\in\Z$, $\ind{T^iL}{K'\cap
  T^iL}=\ind{L}{K'\cap L}$. But then
\begin{multline*}
\ind{K}{K'}=\ind{K}{K\cap K'}\\=\ind{\bigcup_{i=-\infty}^\infty
  T^iL}{K'\cap\bigcup_{i=-\infty}^\infty T^iL}=\ind{L}{K'\cap
  L}<\infty,
\end{multline*}
so we are done.
\end{proof}

\begin{lem}[case of $K/K_c$]\label{find}
Let $K$ be a finitely generated torsion-free abelian group and let $T$
be a hyperbolic automorphism of $K$, so that $K$ has the structure of
a $\lp{t}$-module. Let $A$ be a finite $t$-generating set for $K$ and
let $B=\bigsqcup_{i=-\infty}^\infty t^iA$.  Let $k\in K-\set{0}$. Then
there is a finite subset $C$ of $B$ such that any $t$-word
representing $k$ must contain a letter $\in C$.
\end{lem}

\begin{proof}
Since none of the eigenvalues of $T$ have absolute value $1$,
$K\ten\R$ decomposes as the direct sum of an expanding subspace $K_e$
and a contracting subspace $K_c$. The restrictions of $T$ to $K_c\cap
K$ and of $T\inv $ to $K_e\cap K$ must have all their eigenvalues with
absolute value $<1$. Since $T$ is an automorphism, it follows that
$K_e\cap K=K_c\cap K=\set{0}$.  For every $k\in K$, let $k=k_e+k_c$,
where $k_e\in K_e$ and $k_c\in K_c$.  Let $k\in K$, and suppose at
least one of $k_e$ and $k_c$ is $\in K$.  Then the other must be too,
so $k_e=k_c=0$, so $k=0$.

Choose a norm $\Abs{\cdot}$ on $K\ten\R$ and let $k\in K-\set{0}$. Let
$H_k$ be the Hausdorff distance with respect to this norm between
$k_c$ and $K$.  Since $K$ is finitely generated, hence discrete, and
$k_c\notin K$ by the preceding paragraph, $H_k>0$.

For $w'$ a $t$-word, let $\pi(w')$ be the element of $K$ it
represents.  Let $w$ be a $t$-word of length $l$ in $\pi\inv(k)$,
where $k\in K-\set{0}$.  For any $m\in\Z$, let $w^{>m}$ be the sum of
all letters of $w$ in $\bigsqcup_{i=m+1}^\infty t^iA$ and $w^{<m}$ be
the sum of all letters of $w$ in $\bigsqcup_{-\infty}^{m-1} t^iA$.
Suppose $w=w^{>m}+w^{<-m}$.  Let $k^{>m}=\pi(w^{>m})$ and
$k^{<-m}=\pi(w^{<-m})$.  Then $k^{<-m}=k_c+k^{<-m}_e-k^{>m}_c\in K$.

There are $D\in\R$ and $E>1$ depending only on $K$, $T$ and $A$ such
that, for all $a\in A$, $\Abs{(t^ia)_e}<DE^i$ for all $i\le0$ and
$\Abs{(t^ia)_c}<DE^{-i}$ for all $i\ge0$.  It follows that there are
$F\in\R$ and $E>1$ such that $\Abs{k^{<-m}_e}$ and
$\Abs{k^{>m}_c}<FE^{-m}$ for all $m\in\W$. Putting this together with
the last two paragraphs, we get that $2FE^{-m}>H_k$. We are thus done
if we let $M=-\gint{\log_E[H_k/(2F)]}$ and $C=\bigsqcup_{i=-M}^Mt^iA$.
\end{proof}

\begin{lem}\label{finc}
Let $K$ be a finite-rank torsion-free abelian group and let $L$ be a
full-rank lattice in $K$. Let $T$ be an automorphism of $K$ such that
\begin{itemize}
\item $T$ acts by an endomorphism on $L$,
\item $\bigcup_{i=-\infty}^\infty T^iL=K$ and
\item $\bigcap_{i=-\infty}^\infty T^iL=\set{0}$.
\end{itemize}
Let $A$ be a finite subset of $K$ not including $0$ and let
$B=\bigcup_{i=-\infty}^\infty T^iA$.  Suppose $B$ generates $K$ as a
group and let $\func[\onto]{\pi}{\Z^B}{K}$ be the projection. Let
$\Abs{\cdot}$ denote the $L_1$ norm on $\Z^B$. For $w\in\Z^B$, let
$\Lm(w)$ denote the largest $n\in\Z$ such that every letter of $w$ is
$\in T^nL$.  Then for every $l\in\N$ there is $n$ such that any
nonempty $w\in\pi\inv(T^{\Lm(w)+n}L)$ with $\Abs{w}<l$ has a nonempty
subword $\in\ker\pi$.
\end{lem}

\begin{proof}
The proof is by induction on $l$. The statement is obvious for $l=1$
or $l=2$; let $n=1$.

For $l>2$, suppose $w\in\Z^B$ is nonempty with $\Abs{w}<l$. If
$w\in\ker\pi$, then $w$ is itself a nonempty subword of $w$ in
$\ker\pi$, so assume $w\notin\ker\pi$. For the same reason, assume
none of the letters of $w$ are $0$. Let $\LM(w)$ denote the greatest
$i$ such that some letter of $w$ is in $T^iL$.  Let $b$ be a letter of
$w$ in $T^{\LM(w)}L$.  Let $w'=w-b$, that is the word obtained by
deleting $b$ from $w$.

Suppose $w'\notin\ker\pi$. Then $w'$ is nonempty, so $\Lm(w')=\Lm(w)$.
Also, since $w'$ is a subword of $w$, it has no nonempty subword
$\in\ker\pi$.  Thus, since $\Abs{w'}<l-1$, it must, by induction, be
$\notin\pi\inv(T^{\Lm(w)+n'}L)$, where $n'$ is the $n$ given by
applying the induction assumption. Thus either
$w\notin\pi\inv(T^{\Lm(w)+n'})$ or $b\notin\pi\inv(T^{\Lm(w)+n'})$. In
the latter case, $\LM(w)<\Lm(w)+n'$.

For $i\in\Z$, let
\[
N_{l,i}=\setst{w\in\Z^B-\ker\pi}{\Lm(w)=i,\LM(w)<i+n',\Abs{w}<l}.
\]
Since $A$ is finite, so is $N_{l,i}$.  Let $n''$ be the least integer
such that $w\notin\pi\inv(T^{\Lm(w)+n''}L)$ for all $w\in N_{l,i}$.
Such an integer exists since $N_{l,i}$ is finite and disjoint from
$\ker\pi$.  (It is clear that $n''$ does not depend on $i$.)  Then the
claim follows for $l$, letting $n$ be the greater of $n'$ and $n''$.

We are done by induction.
\end{proof}

\begin{cor}\label{fincc}
Let $K$ be a finite-rank torsion-free abelian group and let $L$ be a
full-rank lattice in $K$. Let $T$ be an automorphism of $K$ such that
\begin{itemize}
\item $T$ acts by an endomorphism on $L$,
\item $\bigcup_{i=-\infty}^\infty T^iL=K$ and
\item $\bigcap_{i=-\infty}^\infty T^iL=\set{0}$.
\end{itemize}
Let $A$ be a finite subset of $K$ and let
$B=\bigcup_{i=-\infty}^\infty T_iA$. Suppose $B$ generates $K$ as a
group and let $\func[\onto]{\pi}{\Z^B}{K}$ be the projection.  Let
$\Abs{\cdot}$ denote the $L_1$ norm on $\Z^B$. Then for every $l\in\N$
and $i\in\Z$ there exist $n_1$ and $n_2\in\Z$ such that each word
$w\in\pi\inv(T^iL-T^{i+1}L)$ with $\Abs{w}<l$ contains a letter in
$T^{n_1}A-T^{n_2}A$.
\end{cor}

\begin{proof}
We may assume $0\notin A$.

Any word $w$ satisfying the conditions contains a subword
\[
w'\in\pi\inv(T^iL-T^{i+1}L)
\]
such that
\begin{itemize}
\item $\Abs{w'}<l$,
\item no nonempty subword of $w'$ is $\in\ker\pi$ and
\item no letter of $w'$ is $\in T^{i+1}L$.
\end{itemize}
By Lemma~\ref{finc}, $w'$ contains all its letters $\in T^{i-n+1}L$.
Since $w'$ is clearly nonempty, we are done.
\end{proof}

\begin{cor}[case of $K/K_d$]\label{finccc}
Let $K$ be a finite-rank torsion-free abelian group and let $L$ be a
full-rank lattice in $K$. Let $T$ be an automorphism of $K$ such that
\begin{itemize}
\item $T$ acts by an endomorphism on $L$,
\item $\bigcup_{i=-\infty}^\infty T^iL=K$ and
\item $\bigcap_{i=-\infty}^\infty T^iL=\set{0}$,
\end{itemize}
so that $K$ has the structure of a $\lp{t}$-module.  Let $A$ be a
finite $t$-generating set for $K$ and let
$B=\bigsqcup_{i=-\infty}^\infty t^iA$.  Let $k\in K-\set{0}$. Then for
every $l\in\N$ there is a finite subset $C$ of $B$ such that any
$t$-word representing $k$ of length $<l$ must contain a letter $\in
C$.
\end{cor}

\begin{proof}
This is just a special case of Corollary~\ref{fincc}.
\end{proof}

We are now ready for the

\begin{proof}[Proof of Proposition~\ref{finite}]
Let $K_d=\bigcap_{i=-\infty}^\infty T^iL$. By Lemma~\ref{gsplit}, $K$
has a finite-index $\lp{t}$-submodule $K'\supseteq K_d$ such that
$K_d$ is a complemented $\lp{t}$-sub\-mod\-ule of $K$. Let $K_c$ be a
complement for $K_d$ in $K'$. Then the actions of $t$ on $K/K_d$ and
$K/K_c$ are hyperbolic. Let $\func[\onto]{\phi_d}{K}{K/K_d}$ and
$\func[\onto]{\phi_c}{K}{K/K_c}$ be the projections.  Let $A_d$ and
$A_c$ denote $\phi_d(A)$ and $\phi_c(A)$ respectively.  Clearly, they
are $t$-generating sets for $K/K_d$ and $K/K_c$ respectively. By a
slight abuse of notation, if $w$ is a $t$-word in $A$, $\phi_d(w)$ and
$\phi(c)(w)$ will denote the corresponding $t$-words in $A_d$ and
$A_c$.

Let $B_c=\bigsqcup_{i=-\infty}^\infty t^iA_c$ and
$B_d=\bigsqcup_{i=-\infty}^\infty t^iA_d$. For $w$ a $t$-word in $A$,
$A_c$ or $A_d$, let $\pi(w)$, $\pi_c(w)$ or $\pi(d)(w)$ be the element
of $K$, $K/K_c$ or $K/K_d$ it represents.

We prove the result by induction on $l$. If $l=1$ then there are only
finitely many $t$-words of length $1$ representing any $k\in
K-\set{0}$, since $T$ is hyperbolic and $A$ is finite. If $l>1$, let
$k\in K-\set{0}$.  Since $K_d\cap K_c=\set{0}$, at least one of
$\phi_d(k)$ and $\phi_c(k)$ is nonzero.  If $w\in\pi\inv(k)$ is a
$t$-word of length $l$ in $A$ then $\phi_d(w)\in\pi_d\inv(\phi_d(k))$
and $\phi_c(w)\in\pi_c\inv(\phi_c(k))$ are $t$-words of length $l$ in
$A_d$ and $A_c$, respectively. But $K/K_c$ is a finite extension of
$K'/K_c\isom K_d$, so it is finitely generated.  Thus, by
Lemma~\ref{find}, if $\phi_c(k)\ne0$ then there is a finite
$C\subset\phi_c(B)$ depending only on $K$, $T$, $A$ and $k$ (not on
$w$) such that $\phi_c(w)$ contains a letter from $C$

If $\phi_c(k)=0$ then $\phi_d(k)\ne0$. Since $K/K_d$ satisfies the
hypotheses of Corollary~\ref{finccc}, there is again some finite
$C\subset\phi_d(B)$ depending only on $K$, $T$, $A$, $k$ and $l$
(again, not on $w$) such that $\phi_d(w)$ contains a letter from $C$.
Putting this together with the preceding paragraph, there is some
finite $C\subset B$ depending only on $K$, $T$, $A$, $k$ and $l$ which
contains a letter from each $t$-word $w\in\pi\inv(k)$ of length at
most $l$.

If $w$ contains a letter $\in\pi\inv(k)$ then the remainder of $w$ is
a nonempty subword representing $0$, so $w$ is not minimal. We may
thus assume that no letter of $w$ represents $k$. In particular, by
the preceding paragraph, $w$ contains a letter $c\in C$ not
representing $k$.  Then $w-c$ is a minimal $t$-word of length at most
$l-1$ representing an element of the finite set $\setst{k-\pi(c)}{c\in
  C}-\set{0}$ with no nonempty subword representing $0$.  We are done
by induction.
\end{proof}

\end{document}